\documentclass[12pt]{iopart}

\usepackage{graphicx}

\usepackage{amsmath}
\usepackage{mathrsfs}
\usepackage{subfigure}
\usepackage{amssymb}

\usepackage{xcolor}

\begin{document}

\title[Solving Large-Scale PDE-constrained Bayesian Inverse Problems with RMHMC]{Solving Large-Scale PDE-constrained Bayesian Inverse Problems with Riemann Manifold Hamiltonian Monte Carlo}


\author{Tan Bui-Thanh$^{\dagger}$
  and Mark Girolami$^{\ddagger}$}

\address{$^{\dagger}$ Department of Aerospace Engineering and Engineering Mechanices \\ Institute for Computational Engineering
    \& Sciences \\ The University of Texas at Austin, Austin, TX 78712,
    USA.}

\address{$^{\ddagger}$ Department of Statistics \\
University of Warwick, Coventry, CV4 7AL, United Kingdom.}
\ead{tanbui@ices.utexas.edu, m.girolami@warwick.ac.uk}

\renewcommand{\thefootnote}{\arabic{footnote}}

\newcommand{\TODO}[1]{ \fbox{\parbox{3in}{\bf TODO: #1}}}

\newcommand{\grbf}[1] {\mbox{\boldmath${#1}$\unboldmath}}
\newcommand{\gbf}[1] {\mathbf{#1}}

\newcommand{\beq} {\begin{equation}}
\newcommand{\eeq} {\end{equation}}
\newcommand{\bdm} {\begin{displaymath}}
\newcommand{\edm} {\end{displaymath}}
\newcommand{\bit}{\begin{itemize}}
\newcommand{\eit}{\end{itemize}}
\newcommand{\bde}{\begin{description}}
\newcommand{\ede}{\end{description}}
\newcommand{\bce}{\begin{center}}
\newcommand{\ece}{\end{center}}
\newcommand{\ben} {\begin{enumerate}}
\newcommand{\een} {\end{enumerate}}
\newcommand{\bea} {\begin{eqnarray}}
\newcommand{\eea} {\end{eqnarray}}
\newcommand{\barr} {\begin{array}}
\newcommand{\earr} {\end{array}}
\newcommand{\bean} {\begin{eqnarray*}}
\newcommand{\eean} {\end{eqnarray*}}
\newcommand{\edoc} {\end{document}}
\newcommand{\On}{\hspace{0.2cm} \mbox{on} \hspace{0.2cm}}
\newcommand{\In}{\hspace{0.2cm} \mbox{in} \hspace{0.2cm}}
\newcommand{\Minimize}{\hspace{0.2cm} \mbox{minimize} \hspace{0.2cm}}
\newcommand{\Maximize}{\hspace{0.2cm} \mbox{maximize} \hspace{0.2cm}}
\newcommand{\SubjectTo}{\hspace{0.2cm} \mbox{subject} \hspace{0.15cm}
            \mbox{to} \hspace{0.2cm}}

\def\etal{{\it et al.~}}

\newcommand{\point}[1] {\ensuremath{\boldsymbol{#1}}}
\newcommand{\vvec}[1] {\ensuremath{\boldsymbol{#1}}}
\renewcommand{\vec}[1] {\ensuremath{\boldsymbol{#1}}}
\newcommand{\cvec}[1] {\ensuremath{\boldsymbol{{#1}}}}
\newcommand{\ten}[1] {\ensuremath{\boldsymbol{#1}}}
\newcommand{\tenfour}[1] {\ensuremath{\boldsymbol{\mathsf{#1}}}}
\newcommand{\comp}[1]{\begin{pmatrix}{ #1 }\end{pmatrix}}
\newcommand{\vv}{\ensuremath{\vec{v}}}
\newcommand{\vr}{\left(\ensuremath{\vec{r}}\right)}
\newcommand{\att}{\left(t\right)}
\newcommand{\tvv}{\ensuremath{\tilde{\vec{v}}}}
\newcommand{\ww}{\ensuremath{\vec{w}}}
\newcommand{\GG}{\ensuremath{\ten{G}}}
\newcommand{\HH}{\ensuremath{\ten{H}}}
\newcommand{\EE}{\ensuremath{\vec{E}}}
\newcommand{\FF}{\ensuremath{\boldsymbol{\mathfrak{F}}}}
\newcommand{\tEE}{\ensuremath{\tilde{\vec{E}}}}
\newcommand{\nn}{\ensuremath{\vec{n}}}
\renewcommand{\SS}{\ensuremath{\ten{S}}}
\newcommand{\tSS}{\ensuremath{\tilde{\ten{S}}}}
\newcommand{\cp}{\ensuremath{{c_p}}}
\newcommand{\cs}{\ensuremath{{c_s}}}
\newcommand{\tcs}{\ensuremath{{\tilde{c}_s}}}
\newcommand{\Vp}{\ensuremath{{V^e_p}}}
\newcommand{\Vs}{\ensuremath{{V^e_s}}}
\newcommand{\Sp}{\ensuremath{{S^e_p}}}
\newcommand{\Ss}{\ensuremath{{S^e_s}}}
\newcommand{\nxnx}[1] {\ensuremath{\nn\times\left(\nn\times{#1}\right)}}

\newcommand{\dd}[2] {\ensuremath{\frac{\partial {#1}}{\partial {#2}}}}
\newcommand{\DD}[2] {\ensuremath{\frac{d {#1}}{d {#2}}}}
\newcommand{\Grad} {\ensuremath{\nabla}}
\newcommand{\Gradv}[1]{\ensuremath{\nabla_{#1}}}
\newcommand{\Div} {\ensuremath{\nabla\cdot}}
\newcommand{\Divv}[1] {\ensuremath{\nabla_{#1}\cdot}}
\newcommand{\Curl} {\ensuremath{\nabla\times}}
\newcommand{\Cten}{\ensuremath{\tenfour{C}}}

\newcommand{\eval}[2][\right]{\relax
  \ifx#1\right\relax \left.\fi#2#1\rvert}

\providecommand{\abs}[1]{\left\lvert#1\right\rvert}
\providecommand{\norm}[1]{\left\lVert#1\right\rVert}

\newcommand{\Nel} {\ensuremath{{N_\text{el}}}}
\newcommand{\Ned} {\ensuremath{{N_\text{ed}}}}
\newcommand{\De} {\ensuremath{{\mathsf{D}^e}}}
\newcommand{\Dep} {\ensuremath{{\mathsf{D}^{e'}}}}
\newcommand{\Dhat} {\ensuremath{\hat{\mathsf{D}}}}
\newcommand{\Dehat} {\ensuremath{{\hat{\mathsf{D}}^e}}}
\newcommand{\half} {\ensuremath{\frac{1}{2}}}
\newcommand{\IN} {\ensuremath{{\mathcal{I}_N}}}
\newcommand{\INe} {\ensuremath{{\mathcal{I}_{N_e}}}}

\newcommand{\mc}[1]{\mathcal{#1}}
\newcommand{\mcC}[1]{\mathcal{#1}}
\newcommand{\mb}[1]{\mathbf{#1}}
\newcommand{\mbb}[1]{\mathbb{#1}}
\newcommand{\mi}[1]{\mathit{#1}}
\newcommand{\mf}[1]{\mathfrak{#1}}
\newcommand{\Oi}[1]{\int_{\mc{D}} #1 \, d\vec{x}}
\newcommand{\TOi}[1]{\int_{0}^T \int_{\mc{D}} #1 \, d\vec{x}dt}
\newcommand{\Ti}[1]{\int_{0}^T #1 \, dt}
\newcommand{\TDei}[1]{\int_{0}^T\int_{\De} #1 \, d\vec{x}dt}
\newcommand{\TDeid}[1]{\int_{0}^T\int_{\De,N} #1 \, d\vec{x}dt}
\newcommand{\Dei}[1]{\int_{\De} #1 \, d\vec{x}}
\newcommand{\DeI}[1]{\int_{D} #1 \, d\vec{x}}
\newcommand{\DeiM}[1]{\int_{\Dhat} #1 \, d\vec{r}}
\newcommand{\DeMd}[1]{\int_{\De,N_e} #1 \, d\vec{x}}
\newcommand{\DeiMd}[1]{\int_{\Dhat,N_e} #1 \, d\vec{r}}
\newcommand{\DeMdN}[1]{\int_{D,N} #1 \, d\vec{x}}
\newcommand{\DeiMdN}[1]{\int_{\Dhat,N} #1 \, d\vec{r}}
\newcommand{\DeiMdNC}[2]{\int_{\Dhat,N_{#2}} #1 \, d\vec{r}}
\newcommand{\pDei}[1]{\int_{\partial \De} #1 \, d\vec{x}}
\newcommand{\pDeiM}[1]{\int_{\partial \Dhat} #1 \, d\vec{r}}
\newcommand{\pDeiMd}[1]{\int_{\partial \Dhat,N_e} #1 \, d\vec{r}}
\newcommand{\pDeiMdN}[1]{\int_{\partial \Dhat,N} #1 \, d\vec{r}}
\newcommand{\pDeiMdNC}[2]{\int_{\partial \Dhat,N_{#2}} #1 \, d\vec{r}}
\newcommand{\pMortar}[2]{\int_{\mc{M}_{#2}} #1 \, d\vec{r}}
\newcommand{\pMortard}[2]{\int_{\mc{M}^{#2},N_{#2}} #1 \, d\vec{r}}
\newcommand{\pMortardd}[3]{\int_{\mc{M}^{#2},N_{#3}} #1 \, d\vec{r}}
\newcommand{\pMortardh}[3]{\int_{\mc{M}_{#2},N_{#3}} #1 \, d\vec{r}}
\newcommand{\pMortardhn}[4]{\int_{{#2}_{#3},N_{#4}} #1 \, d\vec{r}}
\newcommand{\TpDei}[1]{\int_{0}^T\int_{\partial \De} #1 \, d\vec{x}}
\newcommand{\TpDeid}[1]{\int_{0}^T\int_{\partial \De,N} #1 \, d\vec{x}}
\newcommand{\Deid}[1]{\int_{\De,N_e} #1 \, d\vec{x}}
\newcommand{\DeidN}[1]{\int_{\De,N} #1 \, d\vec{x}}
\newcommand{\pDeidN}[1]{\int_{\partial \De,N} #1 \, d\vec{x}}
\newcommand{\pDeid}[1]{\int_{\partial \De,N_e} #1 \, d\vec{x}}
\newcommand{\nor}[1]{\left\| #1 \right\|}
\newcommand{\snor}[1]{\left| #1 \right|}
\newcommand{\LRp}[1]{\left( #1 \right)}
\newcommand{\LRs}[1]{\left[ #1 \right]}
\newcommand{\LRa}[1]{\left< #1 \right>}
\newcommand{\LRc}[1]{\left\{ #1 \right\}}
\newcommand{\pp}[2]{\frac{\partial #1}{\partial #2}}
\newcommand{\ppcp}[1]{\frac{\partial #1}{\partial \vec{c_p}}}
\newcommand{\ppcpj}[1]{\frac{\partial #1}{\partial c_p}}
\newcommand{\ppcpmf}[1]{\frac{\partial #1}{\partial \vec{\mf{c_p}}}}
\newcommand{\ppm}[1]{\frac{\partial #1}{\partial \vec{m}}}
\newcommand{\ppho}[3]{\frac{\partial^{#3} #1}{{\partial #2}^{#3}}}
\newcommand{\ppmi}[1]{\frac{\partial #1}{\partial m_i}}
\newcommand{\ppmj}[1]{\frac{\partial #1}{\partial m_j}}
\newcommand{\ppcpm}[1]{\frac{\partial #1}{\partial \vec{c_p}^-}}
\newcommand{\ppcpp}[1]{\frac{\partial #1}{\partial \vec{c_p}^+}}
\newcommand{\ppcs}[1]{\frac{\partial #1}{\partial \vec{c_s}}}
\newcommand{\ppcsm}[1]{\frac{\partial #1}{\partial \vec{c_s}^-}}
\newcommand{\ppcsp}[1]{\frac{\partial #1}{\partial \vec{c_s}^+}}
\newcommand{\td}[2]{\frac{d #1}{d #2}}
\newcommand{\Rvert}[1]{\left. #1 \right|}
\newcommand{\ipDed}[3]{\LRp{ #1, #2}_{\De,N}^{#3}}

\newcommand{\bigO}{\mathcal{O}}

\newcommand{\iOm}[1]{\int_{\Omega} #1 \, d\Omega}
\newcommand{\iGb}[1]{\int_{\partial\Omega\setminus\Gamma_R} #1 \, ds}
\newcommand{\iGs}[1]{\int_{\Gamma_R} #1 \, ds}
\newcommand{\iGsy}[1]{\int_{\Gamma_{s}} #1 \, ds(\mb{y})}
\newcommand{\RE}{\mbox{{I}}\hspace{-0.5ex}\mbox{{R}}}
\newcommand{\iC}[1]{\int_{0}^{2\pi} #1 \, d\theta}
\newcommand{\iGr}[1]{\int_{\Gamma_{\infty}} #1 \, ds}

\newcommand{\jump}[1] {\ensuremath{\LRs{\![#1]\!}}}
\newcommand{\diff}[1] {\ensuremath{\LRs{#1}}}
\newcommand{\average}[1] {\ensuremath{\LRc{\!\{#1\}\!}}}

\newcounter{tablecell}

\newcommand*{\numbercell}{%
  \stepcounter{tablecell}%
  \thetablecell. 
}

\newtheorem{propo}{Proposition}
\newtheorem{coro}{Corollary}
\newtheorem{theo}{Theorem}
\newtheorem{lem}{Lemma}

\newtheorem{defi}{definition}

\newtheorem{rema}{remark}

\newcommand{\figlab}[1]{\label{fig:#1}}
\newcommand{\eqnlab}[1]{\label{eq:#1}}
\newcommand{\theolab}[1]{\label{theo:#1}}
\newcommand{\corolab}[1]{\label{coro:#1}}
\newcommand{\propolab}[1]{\label{propo:#1}}
\newcommand{\lemlab}[1]{\label{lem:#1}}
\newcommand{\defilab}[1]{\label{defi:#1}}
\newcommand{\remalab}[1]{\label{rema:#1}}
\newcommand{\tablab}[1]{\label{tab:#1}}

\newcommand{\figref}[1]{\ref{fig:#1}}
\newcommand{\theoref}[1]{\ref{theo:#1}}
\newcommand{\defiref}[1]{\ref{defi:#1}}
\newcommand{\remaref}[1]{\ref{rema:#1}}
\newcommand{\cororef}[1]{\ref{coro:#1}}
\newcommand{\proporef}[1]{\ref{propo:#1}}
\newcommand{\lemref}[1]{\ref{lem:#1}}
\newcommand{\eqnref}[1]{\eqref{eq:#1}}
\newcommand{\alglab}[1]{\label{alg:#1}}
\newcommand{\algref}[1]{\ref{alg:#1}}
\newcommand{\seclab}[1]{\label{sect:#1}}
\newcommand{\secref}[1]{\ref{sect:#1}}
\newcommand{\tabref}[1]{\ref{tab:#1}}

\newcommand{\GM}[2]{\mc{N}\left( #1, #2 \right)}

\newcommand{\w}{w}
\newcommand{\wt}{\w^2}
\newcommand{\wth}{\w^3}
\newcommand{\wtth}{\w^{2,3}}
\newcommand{\W}{W}
\renewcommand{\b}{b}
\renewcommand{\u}{u}
\newcommand{\wh}{\hat{\w}}
\newcommand{\uo}{\tilde{\u}}
\newcommand{\ut}{\u^2}
\newcommand{\uth}{\u^3}
\newcommand{\ub}{\mb{\u}}
\renewcommand{\v}{v}
\newcommand{\vb}{\mb{\v}}
\newcommand{\yb}{\mb{y}}
\newcommand{\f}{f}
\newcommand{\g}{g}
\newcommand{\pOmega}{\partial\Omega}
\renewcommand{\d}{d}
\newcommand{\db}{\mb{d}}
\newcommand{\etab}{\boldsymbol{\eta}}
\newcommand{\Ltwo}{L^2\LRp{\Omega}}
\newcommand{\X}{X}
\renewcommand{\L}{{\mb{L}}}
\newcommand{\I}{{\mb{I}}}
\newcommand{\R}{\mathbb{R}}
\newcommand{\A}{\mb{A}}
\newcommand{\B}{B}
\newcommand{\M}{\mb{M}}
\newcommand{\K}{\mb{K}}
\newcommand{\q}{q}
\newcommand{\V}{V}
\newcommand{\zb}{\mb{z}}
\renewcommand{\sb}{\mb{s}}
\newcommand{\rb}{\mb{r}}
\newcommand{\x}{\mb{x}}
\newcommand{\lamh}{\hat{\lambda}}
\newcommand{\lamt}{\lambda^2}
\newcommand{\lamtth}{\lambda^{2,3}}
\newcommand{\lamth}{\lambda^3}
\newcommand{\lamtt}{\tilde{\lamt}}
\newcommand{\like}{{\pi_{\text{like}}\LRp{\db|\u}}}
\newcommand{\post}{{\pi\LRp{\ub | \db}}}
\newcommand{\Ex}{{\mathbb{E}}}
\newcommand{\Q}{\mc{Q}}
\newcommand{\Vb}{\mb{\V}}
\newcommand{\C}{\mb{C}}
\newcommand{\G}{\mb{G}}
\newcommand{\T}{T}
\renewcommand{\H}{\mb{H}}
\renewcommand{\S}{\mb{S}}
\newcommand{\D}{\mb{D}}
\newcommand{\bb}{\mb{b}}
\newcommand{\cb}{\mb{c}}
\newcommand{\J}{\mc{J}}
\newcommand{\Ht}{{\tilde{\H}}}
\newcommand{\Hc}{\mc{H}}
\newcommand{\pb}{\mb{p}}

\begin{abstract}
We consider the Riemann manifold Hamiltonian Monte Carlo (RMHMC) method
for solving statistical inverse problems governed by partial
differential equations (PDEs). The Bayesian framework is employed to
cast the inverse problem into the task of statistical inference whose
solution is the posterior distribution in infinite dimensional
parameter space conditional upon observation data and Gaussian prior measure. We
discretize both the likelihood and the prior using the
$H^1$-conforming finite element method together with a matrix transfer
technique. 
The power of the RMHMC method is that it exploits the
geometric structure induced by the PDE constraints of the underlying inverse
problem. Consequently, each RMHMC posterior sample is almost
uncorrelated/independent from the others providing statistically efficient Markov chain simulation. However this statistical efficiency comes at a computational cost. This motivates us to consider computationally more efficient strategies for
RMHMC. At the heart of our construction is the fact that, Gaussian error structures the Fisher
information matrix coincides with the Gauss-Newton Hessian. We exploit
this fact in considering a computationally simplified RMHMC method combining
state-of-the-art adjoint techniques and the superiority of the RMHMC method. Specifically,
we first form the Gauss-Newton Hessian at the maximum {\it a posteriori}
point and then use it as a fixed constant metric tensor throughout
RMHMC simulation. This eliminates the need for the computationally costly differential geometric Christoffel symbols
which in turn greatly reduces computational effort at a corresponding loss of sampling efficiency. We further reduce the
cost of forming the Fisher information matrix by using a low rank
approximation via a randomized singular value decomposition
technique. This is efficient since a small number of
Hessian-vector products are required. The Hessian-vector product in
turn requires only two extra PDE solves using adjoint
technique. Various numerical results up to $1025$ parameters are
presented to demonstrate the ability of the RMHMC method in exploring
the geometric structure of the problem to propose (almost)
uncorrelated/independent samples that are far away from each other,
and yet the acceptance rate is almost unity. The results also suggest that for the PDE models considered the
proposed fixed metric RMHMC can attain almost as high a quality performance as the original RMHMC, i.e.  generating (almost) uncorrelated/independent samples, while being
two orders of magnitude less computationally expensive.
\end{abstract}



\submitto{\IP}

\maketitle

\section{Introduction}
\seclab{intro}
Inverse problems are ubiquitous in science and engineering. Perhaps
the most popular family of inverse problems is to determine a set of
parameters (or a function) given a set of indirect observations, which
are in turn provided by a parameter-to-observable map plus observation
uncertainties. For example, if one considers the problem of
determining the heat conductivity of a thermal fin given
measured temperature at a few locations on the thermal fin, then: i)
the desired unknown parameter is the distributed heat conductivity,
ii) the observations are the measured temperatures, iii) the
parameter-to-observable map is the mathematical model that describes
the temperature on the thermal fin as a function of the heat
conductivity; indeed the temperature distribution is a solution of an
elliptic partial differential equation (PDE) whose coefficient is the
heat conductivity, and iv) the observation uncertainty is due to the
imperfection of the measurement device and/or model
inadequacy. 

The Bayesian inversion framework refers to a mathematical method that
allows one to solve statistical inverse problems taking into account
all uncertainties in a systematic and coherent manner. The Bayesian
approach does this by reformulating the inverse problem as a problem
in statistical inference, incorporating uncertainties in the
observations, the parameter-to-observable map, and prior information
on the parameter. In particular, we seek a statistical description of
all possible (set of) parameters that conform to the available prior knowledge
and at the same time are consistent with the observations. The
solution of the Bayesian framework is the so-called posterior measure
that encodes the degree of confidence on each set of parameters as the
solution to the inverse problem under consideration.

Mathematically the posterior is a surface in high dimensional
parameter space. The task at hand is therefore to explore the
posterior by, for example, characterizing the mean, the covariance,
and/or higher moments. The nature of this task is to compute high
dimensional integrals for which most contemporary methods are
intractable. Perhaps the most general method to attack these problems
is the Markov chain Monte Carlo (MCMC) method which shall be introduced in subsequent sections.

Let us now summarize the content of the paper. We start with the
description of the statistical inverse problem under consideration in
Section \secref{infiniteBayes}. It is an inverse steady state heat
conduction governed by elliptic PDEs. We postulate a Gaussian measure
prior on the parameter space to ensure that the inverse problem is
well-defined. The prior itself is a well-defined object whose
covariance operator is the inverse of an elliptic differential operator and
with the mean function living in the Cameron-Martin space of the
covariance. The posterior is given by its Radon-Nikodym derivative
with respect to the prior measure, which is proportional to the
likelihood. Since the RMHMC simulation method requires the gradient, Hessian, and the
derivative of the Fisher information operator, we discuss, in some
depth, how to compute the derivatives of the potential function
(the misfit functional) with PDE constraints efficiently using the adjoint
technique in Section \secref{adjoint}. In particular, we define a
Fisher information operator and show that it coincides with the
well-known Gauss-Newton Hessian of the misfit. We next present a
discretization scheme for the infinite Bayesian inverse problem in
Section \secref{FEM}. Specifically, we employ a standard continuous
$H^1$-conforming finite element (FEM) method to discretize both the
likelihood and the Gaussian prior. We choose to numerically compute
the truncated Karhunen-Lo\`eve expansion which requires one to solve an
eigenvalue problem with fractional Laplacian. In order to accomplish
this task, we use a matrix transfer technique (MTT) which leads to a
natural discretization of the Gaussian prior measure.  In Section
\secref{HMC}, we describe the Riemannian manifold Hamiltonian Monte
Carlo (RMHMC) and its variants at length, and its application to our
Bayesian inverse problem. Section \secref{lowrank} presents a low rank
approach to approximate the Fisher information matrix and its inverse
efficiently. This is possibly due to the fact that the Gauss-Newton
Hessian, and hence the Fisher information operator, is a compact
operator. Various numerical results supporting our proposed approach
are presented in Section \secref{results}. We begin this section with
an extensive study and comparison of Riemannian manifold MCMC
methods for problems with two parameters, and end the section with
$1025$-parameter problem. Finally, we conclude the paper in Section
\secref{conclusions} with a discussion on future work.

\section{Problem statement}
\seclab{infiniteBayes} In order to clearly illustrate the challenges arising in PDE-constrained inverse problems for MCMC based Bayesian inference, we consider
the following heat conduction problem governed by an elliptic partial
differential equation in the open and bounded domain $\Omega \subset \R^n$:
\begin{align*}
-\Div\LRp{e^\u\Grad \w} &= 0 & \text{ in } \Omega \\
-e^\u\Grad \w \cdot \mb{n} &= Bi \,u &\text{ on } \partial \Omega \setminus \Gamma_{R}, \\
 -e^\u\Grad \w \cdot \mb{n} &= -1 &\text{ on } \Gamma_{R},
 \end{align*}
where $\w$ is the forward state, $\u$ the logarithm of distributed thermal conductivity on $\Omega$, $\mb{n}$ the unit outward normal on
$\pOmega$, and $Bi$ the Biot number.

In the forward problem, the task is to solve for the temperature
distribution $\w$ given a description of distributed parameter $\u$.
In the inverse problem, the task is to reconstruct $\u$ given some
available observations, e.g, temperature observed at some
parts/locations of the domain $\Omega$. We initially choose to cast
the inverse problem in the framework of PDE-constrained
optimization. To begin, let us consider the following additive
noise-corrupted pointwise observation model\footnote{We assume the
  forward state $\w$ is sufficiently regular, i.e. $\w \in H^{s}, s >
  n/2$, so that $w$ is, by the virtue of the Sobolev embedding
  theorem, continuous, and therefore it is meaningful to measure $w$
  pointwise.}
\begin{equation}\eqnlab{pointwiseObs}
\d_j := \w\LRp{\x_j} + \eta_j, \quad j = 1,\hdots,K,
\end{equation}
where $K$ is the total number of observation locations,
$\LRc{\mb{x}_j}_{j=1}^K$ the set of points at which $\w$ is observed,
$\eta_j$ the additive noise, and $\d_j$ the actual noise-corrupted
observations. In this paper we work with synthetic observations and
hence there is no model inadequacy in
\eqnref{pointwiseObs}. Concatenating all the observations, one can
rewrite \eqnref{pointwiseObs} as
\begin{equation}
\eqnlab{observation}
\db := \mc{G}\LRp{\u} + \etab,
\end{equation}
with $\mc{G} := \LRs{\w\LRp{\mb{x}_1}, \hdots,\w\LRp{\mb{x}_K}}^T$
denoting the map from the distributed parameter $\u$ to the
noise-free observables, $\etab$ being random numbers normally distributed by
$\GM{0}{\L}$ with bounded covariance matrix $\L$, and $\db =
\LRs{\d_1,\hdots,\d_K}^T$. For simplicity, we take $\L = \sigma^2\I$,
where $\I$ is the identity matrix.

Our inverse problem can be now formulated as
\begin{align}
&\min_{\u}  \J\LRp{\u,\db} := \frac{1}{2}\snor{\db-\mc{G}\LRp{\u}}_\L^2 = 
\frac{1}{2\sigma^2} \sum_{j = 1}^K\LRp{\w\LRp{\x_j}-\d_j}^2
\eqnlab{Cost}
\end{align}
\SubjectTo
\begin{subequations}
\eqnlab{forwardEqn}
\begin{align}
-\Div\LRp{e^\u\Grad \w} &= 0 & \text{ in } \Omega, \eqnlab{forward} \\
-e^\u\Grad \w \cdot \mb{n} &= Bi \,u &\text{ on } \partial \Omega \setminus \Gamma_{R}, \eqnlab{forwardbcRobin}\\
 -e^\u\Grad \w \cdot \mb{n} &= -1 &\text{ on } \Gamma_{R},\eqnlab{forwardbcNeumann}.
\end{align}
\end{subequations}
where $\snor{\cdot}_\L:= \snor{\L^{-\half}\cdot}$ denotes the weighted
Euclidean norm induced by the canonical inner product
$\LRp{\cdot,\cdot}$ in $\R^K$. This optimization problem is however
ill-posed. An intuitive reason is that the dimension of observations
$\db$ is much smaller than that of the parameter $\u$, and hence they
provide limited information about the distributed parameter $\u$. As a
result, the null space of the Jacobian of the parameter-to-observation map
$\mc{G}$ is non-empty. Indeed, we have shown that the Gauss-Newton
approximation of the Hessian (which is the square of this Jacobian,
and is also equal to the full Hessian of the data misfit $\mathcal{J}$
evaluated at the optimal parameter) is a compact operator
\cite{Bui-ThanhGhattas12a, Bui-ThanhGhattas12, Bui-ThanhGhattas12f},
and hence its range space is effectively finite-dimensional.

One way to overcome the ill-posedness is to use {\em Tikhonov
  regularization} (see, e.g., \cite{Vogel02}), which proposes to
augment the cost functional \eqnref{Cost} with a quadratic term, i.e.,
\begin{equation}
\tilde{\mc{J}} := \frac{1}{2}\snor{\db-\mc{G}\LRp{\u}}_\L^2 + 
\frac{\kappa}{2}\nor{R^{1/2}\u}^2,
\eqnlab{CostRegularized}
\end{equation}
where $\kappa$ is a regularization parameter, $R$ some regularization
operator, and $\nor{\cdot}$ some appropriate norm.  This method is a
representative of deterministic inverse solution techniques that
typically do not take into account the randomness due to measurements
and other sources, though one can equip the deterministic solution
with a confidence region by post-processing (see, e.g.,
\cite{Vexler04} and references therein). It should be pointed out that
if the regularization term is replaced by the Cameron-Martin norm of
$\u$ (the second term in \eqnref{MAP}), the Tikhonov solution is in
fact identical to the maximum a posteriori point in
\eqnref{MAP}. However, such a point estimate is insufficient for the
purpose of fully taking the randomness into account.

In this paper, we choose to tackle the ill-posedness using a {\em
  Bayesian} framework \cite{Franklin70,
  LehtinenPaivarintaSomersalo89, Lasanen02, Stuart10, Piiroinen05}. We
seek a statistical description of all possible $\u$ that conform to
some prior knowledge and at the same time are consistent with the
observations. The Bayesian approach does this by reformulating the
inverse problem as a problem in {\em statistical inference},
incorporating uncertainties in the observations, the forward map
$\mc{G}$, and prior information. This approach is appealing since it
can incorporate most, if not all, kinds of randomness in a systematic
manner.  To begin, we postulate a Gaussian measure $\mu :=
\GM{\u_0}{\alpha^{-1}\mc{C}}$ on $\u$ in $\Ltwo$ where
\[
\mc{C} := \LRp{I - \Delta}^{-s} =: \mc{A}^{-s}
\]
with the domain of definition
\[
D\LRp{\mc{A}} := \LRc{u \in H^2\LRp{\Omega}: \pp{u}{\mb{n}} = 0 \text{ on } \pOmega},
\]
where $H^2\LRp{\Omega}$ is the usual Sobolev space. Assume that the
mean function $\u_0$ lives in the Cameron-Martin space of $\mc{C}$,
then one can show (see \cite{Stuart10}) that the measure $\mu$ is
well-defined when $s > n/2$ ($d$ is the spatial dimension), and in
that case, any realization from the prior distribution $\mu$ is almost
surely in the H\"older space $\X := C^{0,\beta}\LRp{\Omega}$ with $0 <
\beta < s/2$. That is, $\mu\LRp{X} = 1$, and the
Bayesian posterior measure $\nu$ satisfies the Radon-Nikodym
derivative
\begin{equation}
\eqnlab{RadonNikodym}
\pp{\nu}{\mu}\LRp{\u|\db} \sim \exp\LRp{-\mc{J}\LRp{\u,\db}} = \exp\LRp{-\half\snor{\db -
    \mc{G}\LRp{\u}}^2_\L},
\end{equation}
if $\mc{G}$ is a continuous map from $\X$ to $\R^K$. Note that the
Radon-Nikodym derivative is proportional to the the likelihood defined by
\[
\like \sim \exp\LRp{-\mc{J}\LRp{\u,\db}}.
\]

The maximum a posteriori (MAP) point is
defined as
\begin{equation}
\eqnlab{MAP}
\u^{MAP} := \arg\min_\u \mc{J}\LRp{\u, \db} :=\half\snor{\db -
    \mc{G}\LRp{\u}}^2_\L + \frac{\alpha}{2} \nor{\u}^2_{\mc{C}},
\end{equation}
where $\nor{\cdot}_{\mc{C}} :=
\nor{\mc{C}^{-\half}\cdot}$ denotes the weighted
$\Ltwo$ norm induced by the $\Ltwo$ inner product $\LRa{\cdot,\cdot}$.

\section{Adjoint computation of gradient, Hessian, and the third derivative tensor}
\seclab{adjoint} In this section, we briefly present the adjoint
method to efficiently compute the gradient, Hessian,
and the third derivative of the cost functional \eqnref{Cost}. We
start by considering the weak form of the (first order) forward
equation \eqnref{forwardEqn}:
\begin{equation}
\eqnlab{Forwardw}
\iOm{e^\u\Grad\w \cdot \Grad\lamh} + \iGb{Bi\, \w \lamh} = \iGs{\lamh},
\end{equation}
with $\lamh$ as the test function. Using the standard reduced space
approach (see, e.g., a general discussion in \cite{NocedalWright06}
and a detailed derivation in \cite{Bui-ThanhGhattas13}) one can show
that the gradient $\Grad\mc{J}\LRp{\u}$, namely the Fr\'echet
derivative of the cost functional $\mc{J}$, acting in any direction
$\uo$ is given by
\begin{equation}
\eqnlab{Gradient}
\LRa{\Grad\mc{J}\LRp{\u},\uo} = 
\iOm{\uo e^{\u}\Grad\w\cdot\Grad\lambda},
\end{equation}
where the (first order) adjoint state $\lambda$ satisfies the adjoint equation
\begin{equation}
\eqnlab{Adjointw}
\iOm{e^\u\Grad\lambda \cdot \Grad\wh} + \iGb{Bi\, \lambda \wh} = 
-\frac{1}{\sigma^2}\sum_{j=1}^K\LRp{\w\LRp{\x_j}-\d_j}\wh\LRp{\x_j},
\end{equation}
with $\wh$ as the test function.  On the other hand, the Hessian, the
Fr\'echet derivative of the gradient, acting in directions $\uo$ and
$\ut$ (superscript ``2'' means the second variation direction) reads
\begin{equation}
\eqnlab{FullHessian}
\LRa{\LRa{\Grad^2\mc{J}\LRp{\u},\uo},\ut} =
\iOm{\uo e^\u \Grad\w\cdot\Grad\lamt} + \iOm{\uo e^\u\Grad\wt\cdot\Grad\lambda} + 
\iOm{\uo\ut e^\u\Grad\w\cdot\Grad\lambda},
\end{equation}
where the second order forward state $\wt$ obeys the second order forward equation
\begin{equation}
\eqnlab{Iforward}
\iOm{e^\u\Grad\wt \cdot \Grad\lamh} + \iGb{Bi\, \wt \lamh} = - \iOm{\ut e^\u\Grad\w \cdot \Grad\lamh},
\end{equation}
and the second order  adjoint state $\lamt$ is governed by the second order adjoint equation
\begin{equation}
\eqnlab{Iadjoint}
\iOm{e^\u\Grad\lamt \cdot \Grad\wh} + \iGb{Bi\, \lamt \wh} = 
-\frac{1}{\sigma^2}\sum_{j=1}^K\wt\LRp{\x_j}\wh\LRp{\x_j} - \iOm{\ut e^\u\Grad\lambda \cdot \Grad\wh}.
\end{equation}

We define the generalized Fisher information operator\footnote{Note that the
  Fisher information operator is typically defined for finite dimensional settings in which it is a matrix.}  acting in directions $\uo$ and $\ut$ as
\begin{equation}
\eqnlab{Fisher}
\LRa{\LRa{G\LRp{\u},\uo },\ut} := \Ex_\like\LRs{\LRa{\LRa{\Grad^2\mc{J}\LRp{\u},\uo},\ut}},
\end{equation}
where the expectation is taken with respect to the likelihood---the
distribution of the observation $\db$. Now, substituting
\eqnref{FullHessian} into \eqnref{Fisher} and assuming that the
integrals/derivatives can be interchanged we obtain
\begin{align*}
\LRa{\LRa{G\LRp{\u},\uo },\ut} &= \iOm{\uo e^\u \Grad\w\cdot\Grad\Ex_\like\LRs{\lamt}} + \iOm{\uo e^\u\Grad\wt\cdot\Grad\Ex_\like\LRs{\lambda}} \\
&+ 
\iOm{\uo\ut e^\u\Grad\w\cdot\Grad\Ex_\like\LRs{\lambda}},
\end{align*}
where we have used the assumption that the parameter $\u$ is
independent of observation $\db$ and the fact that $\w$ and $\wt$ do
not depend on $\db$. The next step is to compute
$\Grad\Ex_\like\LRs{\lambda}$ and $\Ex_\like\LRs{\lamt}$. To begin,
let us take the expectation the first order adjoint equation \eqnref{Adjointw} with respect to $\like$ to arrive at
\begin{multline*}
\iOm{e^\u\Grad\Ex_\like\LRs{\lambda} \cdot \Grad\wh} + \iGb{Bi\, \Ex_\like\LRs{\lambda} \wh} = 
\\-\frac{1}{\sigma^2}\sum_{j=1}^K\Ex_\like\LRs{\w\LRp{\x_j}-\d_j}\wh\LRp{\x_j} = 0,
\end{multline*}
where the second equality is obtained from \eqnref{pointwiseObs} and
the assumption $\eta_j \sim \GM{0}{\sigma^2}$. We conclude that 
\begin{equation}
\eqnlab{ExAdjoint}
\Grad\Ex_\like\LRs{\lambda}  = 0.
\end{equation}
On the other hand, if we take the expectation of the second order adjoint equation 
\eqnref{Iadjoint} and use \eqnref{ExAdjoint} we have
\begin{equation}
\eqnlab{ExIadjoint}
\iOm{e^\u\Grad\Ex_\like\LRs{\lamt} \cdot \Grad \wh} + \iGb{Bi\, \Ex_\like\LRs{\lamt} \wh} = 
-\frac{1}{\sigma^2}\sum_{j=1}^K\wt\LRp{\x_j}\wh\LRp{\x_j}.
\end{equation}
Let us define 
\[
\lamtt := \Ex_\like\LRs{\lamt},
\]
then \eqnref{ExIadjoint} becomes
\begin{equation}
\eqnlab{ExIadjointF}
\iOm{e^\u\Grad\lamtt \cdot \Grad\wh} + \iGb{Bi\, \lamtt \wh} = 
-\frac{1}{\sigma^2}\sum_{j=1}^K\wt\LRp{\x_j} \wh\LRp{\x_j}.
\end{equation}

As a result, the Fisher information operator acting along directions
$\uo$ and $\ut$ reads
\begin{equation}
\eqnlab{FisherF}
\LRa{\LRa{G\LRp{\u},\uo },\ut} = \iOm{\uo e^\u \Grad\w\cdot\Grad\lamtt},
\end{equation}
where $\lamtt$ is the solution of \eqnref{ExIadjointF}, a variant of
the second order adjoint equation \eqnref{Iadjoint}. The Fisher
information operator therefore coincides with the Gauss-Newton Hessian of the cost functional \eqnref{Cost}.

The procedure for computing the gradient acting on an arbitrary
direction is clear. One first solves the first order forward equation
\eqnref{Forwardw} for $\w$, then the first order adjoint
\eqnref{Adjointw} for $\lambda$, and finally evaluate
\eqnref{Gradient}. Similarly, one can compute the Hessian (or the
Fisher information operators) acting on two arbitrary directions by
first solving the second order forward equation \eqnref{Iforward} for
$\wt$, then the second order adjoint equation \eqnref{Iadjoint} (or
its variant \eqnref{ExIadjointF}) for $\lamt$ (or $\lamtt$), and
finally evaluating \eqnref{FullHessian} (or \eqnref{FisherF}).

One of the main goals of the paper is to study the Riemann manifold
Hamiltonian Monte Carlo method in the context of Bayesian inverse problems
governed by PDEs. It is therefore essential to compute the derivative
of the Fisher information operator. This task is obvious for problems
with available closed form expressions of the likelihood and the prior,
but it is not so for those governed by PDEs. Nevertheless, using the
adjoint technique we can compute the third order derivative tensor
acting on three arbitrary directions with three extra PDE solves, as
we now show. To that end, recall that the Fisher information operator
acting on directions $\uo$ and $\ut$ is given by \eqnref{FisherF}. The
Fr\'echet derivative of the Fisher information operator along the additional
direction $\uth$ (superscript ``3'' means the third variation direction) is given by
\begin{align}
&\LRa{\LRa{\LRa{\T\LRp{\u}, \uo},\ut}, \uth} := \LRa{\Grad\LRa{\LRa{G\LRp{\u}, \uo},\ut}, \uth} \nonumber \\&= \iOm{\uo\uth e^\u\Grad\w \cdot \Grad\lamtt} + \iOm{\uo e^\u \Grad \wth \cdot \Grad\lamtt} + \iOm{\uo e^\u \Grad\w\cdot\Grad\lamtth}
\eqnlab{tensor},
\end{align}
where $\wth$, $\lamtth$ are the variation of $\w$ and $\lamtt$ in the
direction $\uth$, respectively. One can show that $\uth$ satisfies
another second order forward
equation 
\begin{equation}
\eqnlab{Iforwardn}
\iOm{e^\u\Grad\wth \cdot \Grad \lamh} + \iGb{Bi\, \wth \lamh} = - \iOm{\uth e^\u\Grad\w \cdot \Grad\lamh}.
\end{equation}
Similarly, $\lamtth$ is the solution of the third order adjoint equation
\begin{equation}
\eqnlab{IIadjoint}
\iOm{e^\u\Grad\lamtth \cdot \Grad \wh} + \iGb{Bi\, \lamtth \wh} = 
-\frac{1}{\sigma^2}\sum_{j=1}^K\wtth\LRp{\x_j} \wh\LRp{\x_j}
 - \iOm{\uth e^\u\Grad\lamtt \cdot \Grad\wh},
\end{equation}
and $\wtth$, the variation of $\ut$ in direction $\uth$, satisfies the
following third order forward equation
\begin{align}
&\iOm{e^\u\Grad\wtth \cdot \Grad \lamh} + \iGb{Bi\, \wtth \lamh} = \nonumber\\
&-\iOm{\uth e^\u\Grad\wt \cdot \Grad \lamh}
- \iOm{\uth\ut e^\u\Grad\w \cdot \Grad\lamh}
- \iOm{\ut e^\u\Grad\wth \cdot \Grad\lamh}.
\eqnlab{IIforward}
\end{align}
Note that it would have required four extra PDE solves if one computes
the third derivative of the full Hessian \eqnref{FullHessian}.

It is important to point out that the operator $\T$ is only symmetric
with respect to $\uo$ and $\ut$ since the Fisher information is
symmetric, but not with respect to $\uo$ and $\uth$ or $\ut$ and
$\uth$. The full symmetry only holds for the derivative of the full
Hessian, that is, the true third derivative of the cost functional.

\section{Discretization}
\seclab{FEM} As presented in Section \secref{infiniteBayes}, we view
our inverse problem from an infinite dimensional point of view. As
such, to implement our approach on computers, we need to discretize
the prior, the likelihood and hence the posterior. We choose to use the finite element method. 
In particular, we employ the
standard $H^1\LRp{\Omega}$ finite element method (FEM) to discretize
the forwards and adjoints (the
likelihood), and the operator $\mc{A}$ (the prior). It should be
pointed out that the Cameron-Martin space can be shown (see, e.g.,
\cite{Stuart10}) to be a subspace of the usual fractional Sobolev
space $H^s\LRp{\Omega}$, which is in turn a subspace of
$H^1\LRp{\Omega}$. Thus, we are using a non-conforming FEM approach
(outer approximation). For convenience, we further assume that
the discretized state and parameter live on the same finite
element mesh. Since FEM approximation of elliptic operators is standard (see, e.g., \cite{Ciarlet78}), we will not discuss it here. Instead, we describe the matrix transfer technique (see, e.g, \cite{IlicLiuTurnerEtAl05} and the references therein) to
discretize the prior.

Define $\Q := \mc{C}^{1/2} = \mc{A}^{-s/2}$, then the eigenpairs
$\LRp{\lambda_i,\v_i}$ of $\Q$ define the Karhunen-Lo\`eve (KL)
  expansion of the prior distribution as
\[
\u = \u_0+ \frac{1}{\sqrt{\alpha}}\sum_{i=1}^\infty a_i \lambda_i \v_i, 
\]
where $a_i \sim \GM{0}{1}$. We need to solve
\[
\Q \v_i = \lambda_i\v_i,
\]
or equivalently
\begin{equation}
\eqnlab{eigenProblem}
\mc{A}^{s/2} \v_i = \frac{1}{\lambda_i}\v_i.
\end{equation}

To solve \eqnref{eigenProblem} using the matrix transfer technique (MTT),
let us denote by $\M$ the mass matrix, and $\K$ the stiffness matrix
resulting from the discretization of the Laplacian $\Delta$. The
representation of $\mc{A}$ in the finite element space (see, e.g., \cite{Simpson08} and the references therein)  is given by
\[
\A := \M^{-1}\K + \I.
\]
Let bold symbols denote the corresponding vector of FEM nodal values,
e.g., $\ub$ is the vector containing all FEM nodal values of
$\u$. If we define $\LRp{\sigma_i,\vb_i}$ as
eigenpairs for $\A$, i.e,
\begin{equation}
\eqnlab{eigenmatrix}
\A \vb_i = \sigma_i \vb_i, \text{ or } \A \Vb = \boldsymbol{\Sigma}\Vb
\end{equation}
where $\vb_i^T\M\vb_j = \delta_{ij}$, and hence $\Vb^{-1} = \Vb^T\M$,
$\delta_{ij}$ is the Kronecker delta function, and $\boldsymbol{\Sigma}$ is the
diagonal matrix with entries $\sigma_i$. Since $\A$ is similar to $\M^{-\half}\LRp{\K +
  \M}\M^{-\half}$, a symmetric positive definite matrix, $\A$ has
positive eigenvalues.  Using MTT method, the matrix representation of
\eqnref{eigenProblem} reads
\[
\A^{s/2}\vb_i = \frac{1}{\lambda_i}\vb_i,
\]
where 
\[
\A^{s/2} := \Vb\boldsymbol{\Sigma}^{s/2}\Vb^{-1}.
\]
It follows that 
\[
\lambda_i = \sigma^{-s/2}_i.
\]
The Galerkin FEM approximation of the prior via truncated KL expansion reads
\begin{equation}
\eqnlab{uKLtruncated}
\ub = \ub_0+ \frac{1}{\sqrt{\alpha}}\sum_{i=1}^{N}a_i\lambda_i\vb_i,
\end{equation}
with $\ub$ as the FEM nodal value of the approximate prior sample $\u$
and $N$ as the number of FEM nodal points. Note that for ease in
writing, we have used the same notation $\u$ for both infinite
dimensional prior sample and its FEM approximation. Since $\u
\in \Ltwo$, $\ub$ naturally lives in $\R^N_\M$, the Euclidean
space with weighted inner product $\LRp{\cdot,\cdot}_\M := \LRp{\cdot,\M\cdot}$.

A question arises: what is the distribution of $\ub$? Clearly
$\ub$ is a Gaussian with mean $\ub_0$ since $a_i$ are. The
covariance matrix $\C$ for $\ub$ is defined by
\[
\LRp{\zb, \C\yb}_\M := \Ex\LRs{\LRp{\ub-\ub_0,\M\zb}\LRp{\ub-\ub_0,\M\yb}} = \frac{1}{\alpha}\zb^T \M \Vb \boldsymbol{\Lambda}^2 \Vb^T\M\yb,
\]
where we have used \eqnref{uKLtruncated} to obtain the second equality and  $\boldsymbol{\Lambda}$ is the diagonal matrix with entries $\Lambda_{ii}
= \lambda_i^{-1}$. It follows that 
\begin{equation}
\eqnlab{eqnC}
\C = \frac{1}{\alpha} \Vb \boldsymbol{\Lambda}^2 \Vb^T\M
\end{equation}
as a map from $\R^N_\M$ to $\R^N_\M$, and its inverse can be shown to be
\[
\C^{-1} = \alpha \Vb \boldsymbol{\Lambda}^{-2} \Vb^T\M,
\]
whence the distribution of $\ub$ is
\begin{equation}
\eqnlab{uDist}
\ub \sim \GM{\ub_0}{\alpha \Vb \boldsymbol{\Lambda}^{-2} \Vb^T\M} \sim
\exp\LRs{-\frac{\alpha}{2} \LRp{\ub-\ub_0}^T\M\Vb \boldsymbol{\Lambda}^{-2} \Vb^T\M\LRp{\ub-\ub_0}}.
\end{equation}
As a result, the FEM discretization of the prior can be written as
\begin{align*} 
\frac{\alpha}{2}\nor{\u-\u_0}^2_{\mc{C}} &:= \frac{\alpha}{2} \nor{\mc{A}^{s/2}\LRp{\u-\u_0}}^2 \stackrel{\text{MTT}}{\approx} 
 \frac{\alpha}{2}\LRp{\ub-\ub_0}^T\M\Vb\boldsymbol{\Lambda}^{-2}\Vb^T\M\LRp{\ub-\ub_0}.
\end{align*}
Thus, the FEM approximation of the posterior is given by
\[
\post \sim \exp\LRp{-\half\snor{\db -
    \mc{G}\LRp{\u}}^2_\L} \times \exp\LRp{-\frac{\alpha}{2}\LRp{\ub-\ub_0}^T\M\Vb\boldsymbol{\Lambda}^{-2}\Vb^T\M\LRp{\ub-\ub_0}}.
\]
The detailed derivation of the FEM approximation of infinite Bayesian
inverse problems in general and the prior in particular will be
presented elsewhere \cite{Bui-Thanh14}.



\section{Riemannian manifold Langevin and Hamiltonian Monte Carlo methods}
\seclab{HMC} In this section we give a brief overview of the MCMC
algorithms that we consider in this work. Some familiarity with the
concepts of MCMC is required by the reader since an introduction to
the subject is out of the scope of this paper.

\subsection{Metropolis-Hastings}
For a random vector $\ub \in \mathbb{R}^N$ with
density $\pi(\ub)$ the Metropolis-Hastings algorithm
employs a proposal mechanism
$q(\ub^{*}|\ub^{t-1})$ and proposed
moves are accepted with probability $\min
\left\{1,\pi(\ub^{*})
q(\ub^{t-1}|\ub^{*})/\pi(\ub^{t-1})q(\ub^{*}|\ub^{t-1})\right\}$. Tuning
the Metropolis-Hastings algorithm involves selecting an appropriate
proposal mechanism. A common choice is to use a Gaussian proposal of
the form $q(\ub^{*}|\ub^{t-1}) =
\mathcal{N}(\ub^{*}|\ub^{t-1},\boldsymbol{\Sigma})$,
where $\mathcal{N}(\cdot|\boldsymbol{\mu},\boldsymbol{\Sigma})$
denotes the multivariate normal density with mean $\boldsymbol{\mu}$
and covariance matrix $\boldsymbol{\Sigma}$.

Selecting the covariance matrix however, is far from trivial in most
of cases since knowledge about the target density is
required. Therefore a more simplified proposal mechanism is often
considered where the covariance matrix is replaced with a diagonal
matrix such as $\boldsymbol{\Sigma}=\epsilon\I$ where the
value of the scale parameter $\epsilon$ has to be tuned in order to
achieve fast convergence and good mixing. Small values of $\epsilon$
imply small transitions and result in high acceptance rates while the
mixing of the Markov Chain is poor. Large values on the other hand,
allow for large transitions but they result in most of the samples
being rejected. Tuning the scale parameter becomes even more difficult
in problems where the standard deviations of the marginal posteriors
differ substantially, since different scales are required for each
dimension, and when correlations between different variables exist. In
the case of PDE-constrained inverse problems in very high dimensions
with strong nonlinear interactions inducing complex non-convex
structures in the target posterior this tuning procedure is typically
doomed to failure of convergence and mixing.

There have been many subsequent developments of this basic algorithm
however the most important with regard to inverse problems, arguably,
is the formal definition of Metropolis Hastings in an infinite
dimensional functional space. One of the main failings of Metropolis
Hastings is the drop-off in acceptance probability as the dimension of
the problem increases. By defining the Metropolis acceptance
probability in the appropriate Hilbert space the acceptance
probability should then be invariant to the dimension of the problem
and this is indeed the case as is described in a number of scenarios
by \cite{CotterRobertsStuartEtAl13}. Furthermore the definition of a Markov chain
transition kernel directly in the Hilbert space which exploits
Hamiltonian dynamics in the proposal mechanism followed in
\cite{OttobrePillaiPinskiEtAl14}.

These are important methodological advances for MCMC applied to
Inverse Problems. As the infinite dimensional nature of the problem is
a fundamental aspect of the problem it is sensible that this
characteristic is embedded in the MCMC scheme. In a similar vein by
noting that the statistical model associated with the specific inverse
problem is generated from an underlying partial differential equation
or system of ordinary differential equations a natural geometric
structure structure on the space of probability distributions is
induced. This structure provides a rich source of model specific
information that can be exploited in devising MCMC schemes that are
informed by the underlying structure of the model itself.

In \cite{GirolamiCalderhead11} a way around this situation was
provided by accepting that the statistical model can itself be
considered as an object with an underlying geometric structure that
could be embedded into the proposal mechanism. A class of MCMC methods
were developed based on the differential geometric concepts underlying
Riemannian manifolds.

\subsection{Riemann Manifold Metropolis Adjusted Langevin Algorithm}
Denoting the log of the target density as $\mathcal{L}(\ub) = \log \pi(\ub)$, the manifold Metropolis Adjusted Langevin Algorithm (mMALA) method, \cite{GirolamiCalderhead11}, defines a Langevin diffusion with stationary distribution $\pi(\ub)$ on the Riemann manifold of density functions with metric tensor $\G(\ub)$. By employing a first order Euler integrator for discretising the stochastic differential equation a proposal mechanism with density 
$q(\ub^*|\ub^{t-1}) = \mathcal{N}(\ub^*| \boldsymbol{\mu}(\ub^{t-1},\epsilon),\epsilon^2\G^{-1}(\ub^{t-1}))$
is defined, where $\epsilon$ is the integration step size, a parameter which needs to be tuned, and the $k$th component of the mean function $\boldsymbol{\mu}(\ub,\epsilon)_k$ is
\begin{eqnarray}
\boldsymbol{\mu}(\ub,\epsilon)_k & = & \ub_k +
\frac{\epsilon^2}{2}\left(\G^{-1}(\ub)\nabla_{\ub}\mathcal{L}(\ub)\right)_k
- \epsilon^2
\sum_{i=1}^N\sum_{j=1}^N\G(\ub)_{i,j}^{-1}\Gamma_{i,j}^k \label{eq:meanmMALA}\end{eqnarray}
where $\Gamma_{i,j}^k$ are the Christoffel symbols of the metric in
local coordinates. Note that we have used the Christoffel symbols to
express the derivatives of the metric tensor, and they are computed using the adjoint method presented in Section \secref{adjoint}.

Due to the discretisation error introduced by the first order approximation convergence to the stationary distribution is not guaranteed anymore and thus the Metropolis-Hastings ratio is employed to correct for this bias. In  \cite{GirolamiCalderhead11} a number of examples are provided illustrating the potential of such a scheme for challenging inference problems.

One can interpret the proposal mechanism of RMMALA as a local Gaussian approximation to the target density where the effective covariance matrix in RMMALA is the inverse of the metric tensor evaluated at the current position. Furthermore a simplified version of the RMMALA algorithm, termed sRMMALA, can also be derived by assuming a manifold with constant curvature thus cancelling the last term in Equation (\ref{eq:meanmMALA}) which depends on the Christoffel symbols. Whilst this is a step forward in that much information about the target density is now embedded in the proposal mechanism it is still driven by a random walk. The next approach to be taken goes beyond the direct and scaled random walk by defining proposals which follow the geodesic flows on the manifold of densities and thus presents a potentially really powerful scheme to explore posterior distributions.

\subsection{Riemann Manifold Hamiltonian Monte Carlo}
The Riemann manifold Hamiltonian Monte Carlo (RMHMC) method defines a Hamiltonian on the Riemann manifold of probability density functions by introducing the auxiliary variables $\pb\sim \mathcal{N}(\boldsymbol{0},\G(\ub))$ which are interpreted as the momentum at a particular position $\ub$ and by considering the negative log of the target density as a potential function. More formally the Hamiltonian defined on the Riemann manifold is
\begin{equation}
\Hc(\ub,\pb) = -\mathcal{L}(\ub) +\frac{1}{2}\log\left(2\pi|\G(\ub)|\right) + \frac{1}{2}\pb^T\G(\ub)^{-1}\pb
\end{equation}
where the terms $-\mathcal{L}(\ub) +\frac{1}{2}\log\left(2\pi|\G(\ub)|\right) $ and $\frac{1}{2}\pb^T\G(\ub)^{-1}\pb$ are the potential energy and kinetic energy terms respectively.
 and the dynamics given by Hamiltons equations are
\begin{eqnarray}
\frac{d \ub_k}{dt} &= \frac{\partial \Hc }{\partial \pb_k} = \left(\G(\ub)^{-1}\pb\right)_k \\
\frac{d\pb_k}{dt} &= -\frac{\partial \Hc }{\partial \ub_k} = \frac{\partial \mathcal{L}(\ub)}{\partial \ub_k} -\frac{1}{2}Tr\left[\G(\ub)^{-1}\frac{\partial \G(\ub)}{\partial \ub_k} \right] \nonumber\\
&+\frac{1}{2}\pb^T\G(\ub)^{-1}\frac{\partial \G(\ub)}{\partial \ub_k}\G(\ub)^{-1}\pb 
\eqnlab{Hamilton}
\end{eqnarray}

These dynamics define geodesic flows at a particular energy level and as such make proposals which follow deterministically the most efficient path across the manifold from the current density to the proposed one. Simulating the Hamiltonian requires a time-reversible and volume preserving numerical integrator. For this purpose the Generalised Leapfrog algorithm can be employed and provides a deterministic proposal mechanism for simulating from the conditional distribution, i.e. $\ub^*|\pb \sim \pi(\ub^*|\pb)$. More details about the Generalised Leapfrog integrator can be found in  \cite{GirolamiCalderhead11}. To simulate a path (which turns out to be a local geodesic) across the manifold, the Leapfrog integrator is iterated $L$ times which along with the integration step size $\epsilon$ are parameters requiring tuning. Again due to the discrete integration errors on simulating the Hamiltonian in order to ensure convergence to the stationary distribution the Metropolis-Hastings acceptance ratio is applied. 

The RMHMC method has been shown to be highly effective in sampling from posteriors induced by complex statistical models and offers the means to efficiently explore the hugely complex and high dimensional posteriors associated with PDE-constrained inverse problems.

\section{Low rank approximation of the Fisher information matrix}
\seclab{lowrank} As presented in Section \secref{HMC}, we use the
Fisher information matrix at the MAP point augmented with the Hessian
of the prior as the metric tensor in our HMC simulations. It is
therefore necessary to compute the augmented Fisher matrix and its
inverse. In~\cite{Bui-ThanhGhattas12a, Bui-ThanhGhattas12,
  Bui-ThanhGhattas12f}, we have shown that the Gauss-Newton Hessian of
the cost functional \eqnref{Cost}, also known as the data misfit, is a
compact operator, and that for smooth $\u$ its eigenvalues decay
exponentially to zero. Thus, the range space of the Gauss-Newton
Hessian is effectively finite-dimensional even before discretization,
i.e., it is independent of the mesh. In other words, the Fisher
information matrix admits accurate low rank approximations and the accuracy can be improved as desired by simply increasing the rank of the approximation. We shall
exploit this fact to compute the augmented Fisher information matrix
and its inverse efficiently. We start with the augmented Fisher
information matrix in $\R^N_\M$
\begin{align*}
\G := \M^{-1}\H + \alpha \Vb \boldsymbol{\Lambda}^{-2} \Vb^T\M = \alpha
\Vb\boldsymbol{\Lambda}^{-1}\LRp{\frac{1}{\alpha}\boldsymbol{\Lambda}\Vb^T\H\Vb\boldsymbol{\Lambda}
  + \I}\boldsymbol{\Lambda}^{-1}\Vb^{-1},
\end{align*}
where $\H$ is the Fisher information matrix obtained from
\eqnref{FisherF} by taking $\uo$ and $\ut$ as FEM basis functions.

Assume that $\H$ is compact (see, e.g., \cite{Bui-ThanhGhattas12a,
  Bui-ThanhGhattas12}), together with the fact that
$\boldsymbol{\Lambda}_{ii}$ decays to zero, we conclude that the
prior-preconditioned Fisher information matrix
\[
\Ht := \frac{1}{\alpha}\boldsymbol{\Lambda}\Vb^T\H\Vb\boldsymbol{\Lambda}
\] also has
eigenvalues decaying to zero. Therefore it is expected that the
eigenvalues of the prior-preconditioned matrix decays faster than
those of the original matrix $\H$. Indeed, the numerical results in
Section \secref{results} will confirm this observation.  It follows that (see,
e.g., \cite{FlathWilcoxAkcelikEtAl11, Bui-ThanhGhattasMartinEtAl13}
for similar decomposition) $\Ht$ admits a $r$-rank approximation of the
form
\[
\Ht = \frac{1}{\alpha}\boldsymbol{\Lambda}\Vb^T\H\Vb\boldsymbol{\Lambda} \approx
\Vb_r\S\Vb_r^T, 
\]
where $\Vb_r$ and $\S$ (diagonal matrix) contain the first $r$
dominant eigenvectors and eigenvalues of $\Ht$, respectively. In this
work, similar to \cite{Bui-ThanhGhattasMartinEtAl13}, we use the
one-pass randomized algorithm in \cite{HalkoMartinssonTropp11} to
compute the low rank approximation.  Consequently, the augmented Fisher information matrix becomes
\[
\G \approx \alpha \Vb\boldsymbol{\Lambda}^{-1}\LRp{\Vb_r\S\Vb_r^T
  + \I}\boldsymbol{\Lambda}^{-1}\Vb^{-1}, 
\]
from which we obtain the inverse, by using the Woodbury formula \cite{GolubVan96},
\[
\G^{-1} \approx \frac{1}{\alpha} \Vb\boldsymbol{\Lambda}\LRp{\I - \Vb_r\D\Vb_r^T
  }\boldsymbol{\Lambda}\Vb^{-1}, 
\]
where $\D$ is a diagonal matrix with  $\D_{ii} = \S_{ii}/\LRp{\S_{ii}+1}$.

In the RMHMC method, we need to randomly draw the momentum variable as $\pb \sim
\GM{\mb{0}}{\G}$. If one considers
\[
\pb = \sqrt{\alpha} \Vb \Lambda^{-1}\bb + \sqrt{\alpha} \Vb\boldsymbol{\Lambda}^{-1}\Vb_r\S^{1/2}\cb,
\]
where $\bb_i,\cb_i \sim \GM{0}{1}$, then one can show , by inspection, that
$\pb$ is distributed by $\GM{\mb{0}}{\G}$.

\section{Numerical results}
\seclab{results} 

For convenience, let us recall that the finite element (FEM) approximation of the
posterior is given as
\begin{equation}
\eqnlab{finalPosterior}
\post \sim \exp\LRp{-\half\snor{\db -
    \mc{G}\LRp{\u}}^2_\L} \times \exp\LRp{-\frac{\alpha}{2}\LRp{\ub-\ub_0}^T\M\Vb\boldsymbol{\Lambda}^2\Vb^T\M\LRp{\ub-\ub_0}}, 
\end{equation}
$\ub_0$ is the FEM nodal value of the prior mean function $\u_0$, $\M$
is the mass matrix, $\Vb$ the matrix of eigenvectors defined in
\eqnref{eigenmatrix}, $\boldsymbol{\Lambda}$ the diagonal matrix introduced in
\eqnref{eqnC}, $\L = \sigma^2\I$, $\db$ vector of observation data, and $\mc{G}\LRp{\u}$
the forward map given by the forward equation
\begin{align*}
-\Div\LRp{e^\u\Grad \w} &= 0 & \text{ in } \Omega \\
-e^\u\Grad \w \cdot \mb{n} &= Bi \,u &\text{ on } \partial \Omega \setminus \Gamma_{R}, \\
 -e^\u\Grad \w \cdot \mb{n} &= -1 &\text{ on } \Gamma_{R},
 \end{align*}
that is discretized by the $H^1$-conforming FEM method.

In this section, we study Riemann manifold Monte Carlo methods and
their variations to explore the posterior \eqnref{finalPosterior}. In
particular, we compare the performance of four methods: i) sRMMALA
obtained by ignoring the third derivative in RMMALA, ii) RMMALA, iii)
sRMHMC obtained by first computing the augmented Fisher metric tensor
at the MAP point and then using it as the constant metric tensor, iv)
RMHMC. For all methods, we form the augmented Fisher information
matrix exactly using \eqnref{FisherF} with $\uo$ and $\ut$ as finite
element basis vectors. For RMMALA and RMHMC we also need the
derivative of the metric tensor which is a third order tensor.  It can
be constructed exactly using \eqnref{tensor} with $\uo, \ut$ and
$\uth$ as finite element basis vectors. We also need extra work for
the RMHMC method since each Stormer-Verlet step requires an implicit
solve for both the first half of momentum and full position. For
inverse problems such as those considered in this paper, the fixed
point approach proposed in \cite{GirolamiCalderhead11} does not seem
to converge. We therefore have to resort to a full Newton
method. Since we explicitly construct the metric tensor and its
derivative, it is straightforward for us to develop the Newton
scheme. For all problems considered in this section, we have observed
that it takes at most five Newton iterations to converge.

Note that we limit ourselves in comparing these four methods in the
Riemannian manifold MCMC sampling family. Clearly, other methods are available,
we avoid ``unmatched comparison'' in terms of cost and the level of
exploiting the structure of the problem. Even in this limited family,
RMHMC is most expensive since it  requires not only third derivatives
but also implicit solves, but the ability in generating almost
independent samples is attractive and appealing as we shall show.

Though our proposed approach described in previous sections are valid
for any spatial dimension $d$, we restrict ourselves to a one
dimensional problem, i.e. $d = 1$, to clearly illustrate our points
and findings. In particular, we take $\Omega = \LRs{0,1}$, $\Gamma_R =
\LRc{1}$. We set $Bi = 0.1$ for all examples. As discussed in Section \secref{infiniteBayes},
for the Gaussian prior to be well-defined, we take $s = 0.6 > n/2 = 1/2$. 

\subsection{Two-parameter examples}
\seclab{twoparameter}

We start our numerical experiments with two parameters. This will help
demonstrate various aspects of RMHMC which are otherwise too
computationally expensive for high dimensional problems. In
particular, two-parameter example allows us to compute the complete
third derivative tensor and perform the Newton method for each
Stormer-Verlet step. This in turn allows us to show the capability of the full
RMHMC over its simplified variants in tackling challenging posterior
densities in which the metric tensor changes rapidly.

In order to construct the case with two parameters we consider FEM with 
one finite element. We assume that there is one observation point,
i.e. $K = 1$, and it is placed at the left boundary $x = 0$. In the
first example, we first take $s = 0.6$, $\sigma = 0.1$, and $\alpha =
0.1$. The posterior in this case is shown in Figure
\figref{posterior1}. We start by taking a time step of $0.02$ with $100$
Stormer-Verlet steps for both sRMHMC and RMHMC. The acceptance rate
for both methods is $1$. One would take a time step of $2$ for both
sRMMALA and RMMALA to be comparable with sRMHMC and RMHMC, but the
acceptance would be zero. Instead we take time step of $1$ so that the
acceptance rate is about $0.5$ for sRMMALA and $0.3$ for RMMALA. The
MAP point is chosen as the initial state for all the chains with
$5000$ sample excluding the first $100$ burn-ins. The result is shown
in Figure \figref{comparison1}.

\begin{figure}[h!tb]
  \subfigure[$s = 0.6$, $\alpha = 0.1$, and $\sigma = 0.1$]{
    \includegraphics[trim=1cm 6.5cm 2cm 7.3cm,clip=true,width=0.32\columnwidth]{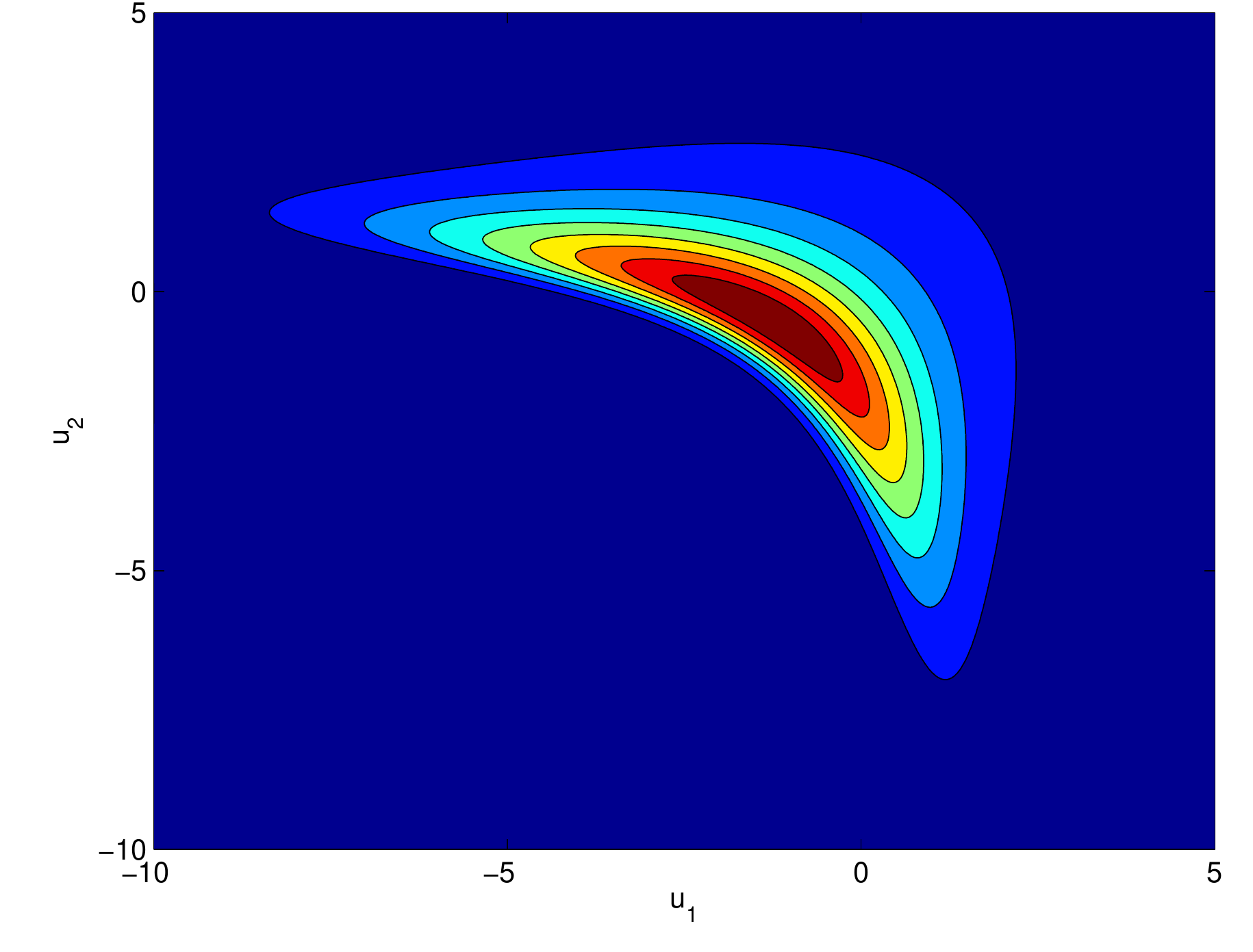}
    \figlab{posterior1}
  }
  \subfigure[$s =
    0.6$, $\alpha = 1$, and $\sigma = 0.01$]{
    \includegraphics[trim=1cm 6.5cm 2cm 7.3cm,clip=true,width=0.32\columnwidth]{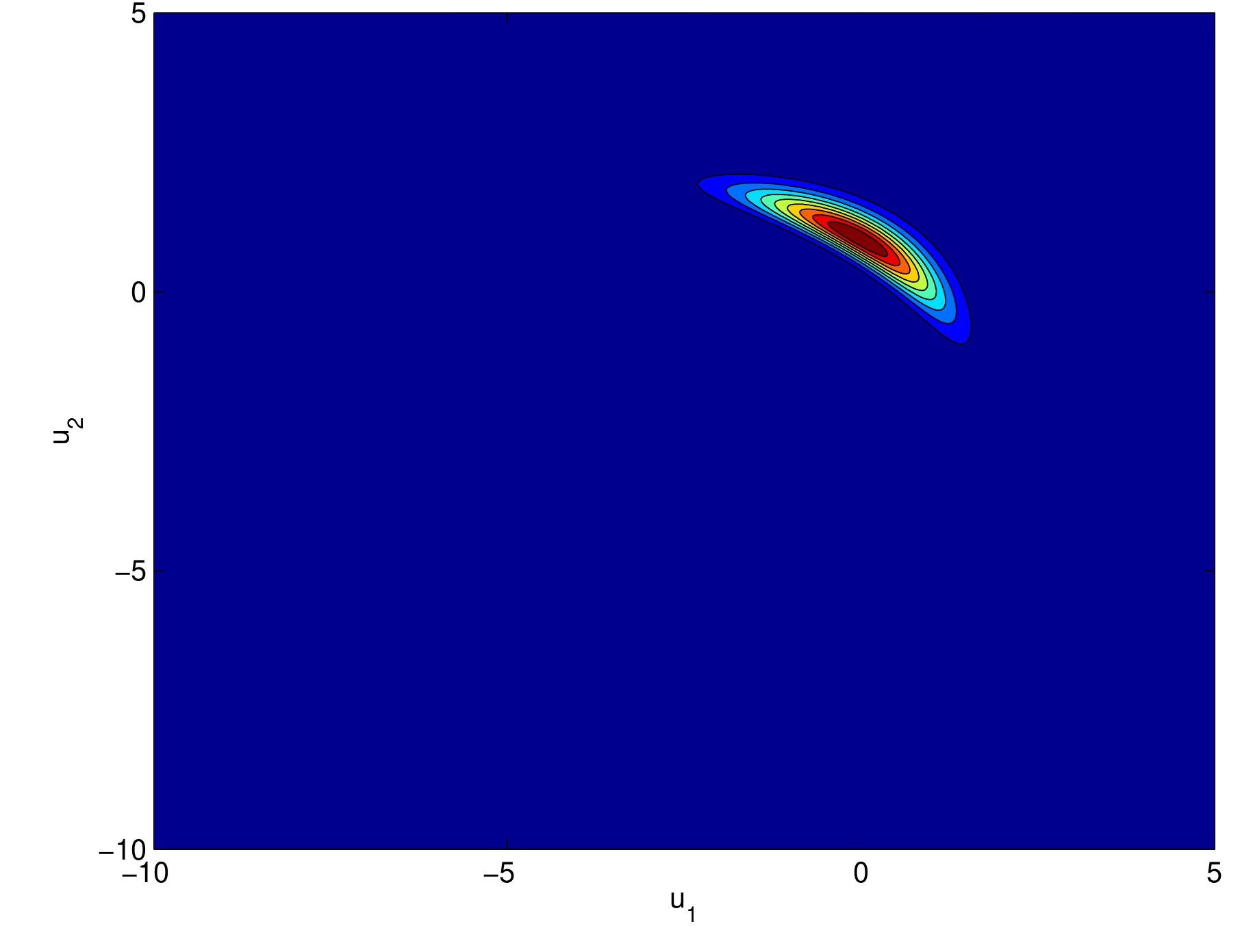}
    \figlab{posterior2}
  }
  \subfigure[$s =
    0.6$, $\alpha = 0.1$, and $\sigma = 0.01$]{
    \includegraphics[trim=1cm 6.5cm 2cm 7.3cm,clip=true,width=0.32\columnwidth]{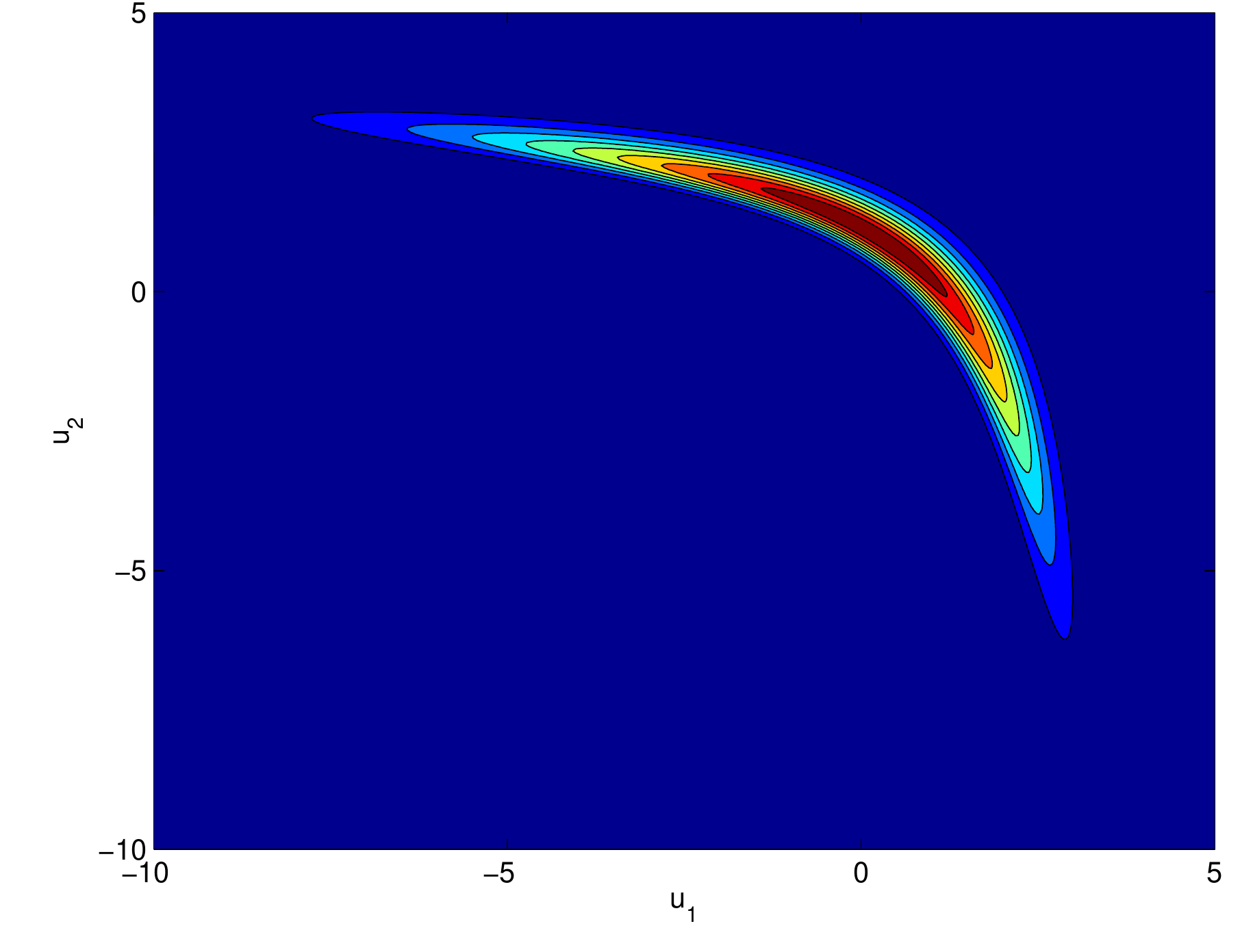}
    \figlab{posterior3}
}
  \caption{The contours of the posterior for three combinations of $s$
    (the prior smoothness), $\alpha$ (the ``amount'' of the prior), and $\sigma$ (the noise standard deviation).} \figlab{Posterior1}
\end{figure}

\begin{figure}[h!tb]
    \includegraphics[trim=1cm 6.5cm 2cm 6.5cm,clip=true,width=0.97\columnwidth]{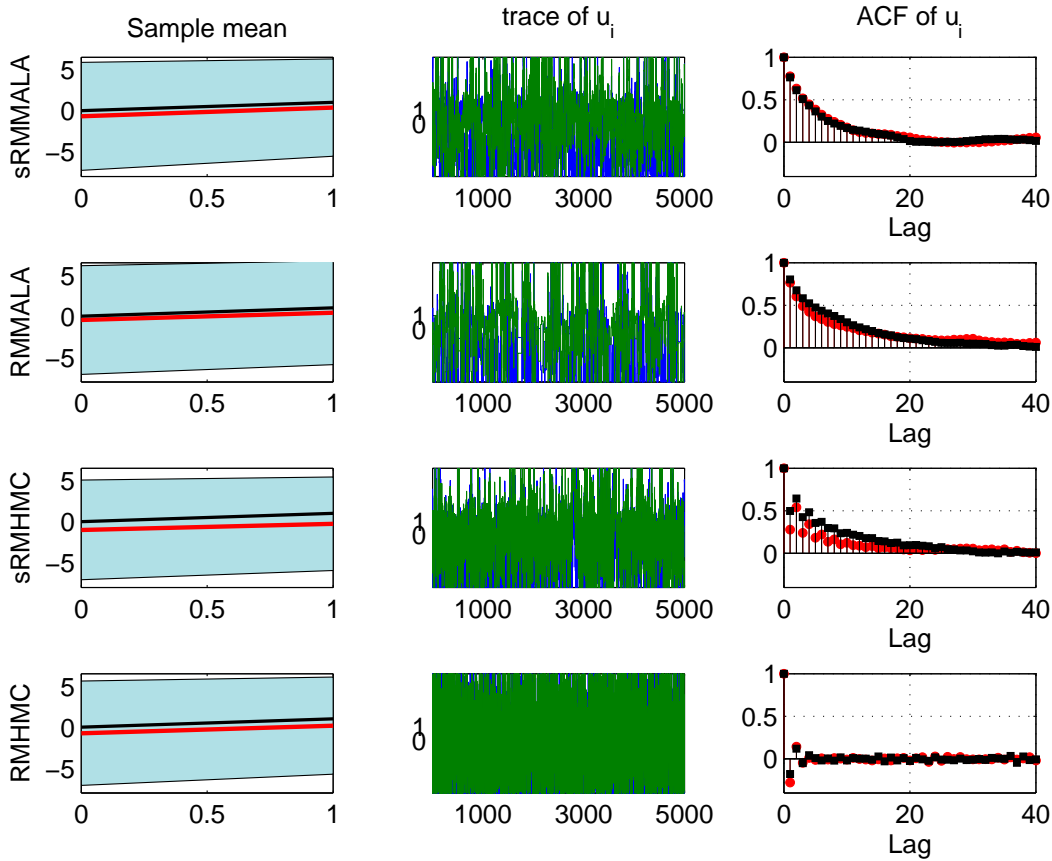}
  \caption{Comparison of simRMMALA, RMMALA, simRMHMC, and RMHMC:
    chains with 5000 samples, burn-in of $100$, starting at the MAP
    point. In this example, $s = 0.6$, $\alpha = 0.1$, and $\sigma =
    0.1$. Time step is $\varepsilon = 1$ for simRMMALA and RMMALA, and
    $\varepsilon = 0.02$ with the number of time steps $L = 100$ for
    simRMHMC and RMHMC. In the left column: the exact synthetic
    solution is black, the sample mean is red, and the shaded region
    is the $95\%$ credibility region. In the middle column: blue is
    the trace plot for $\u_1$ while green is for $\u_2$. In the right
    column: red and black are the autocorrelation function for $\u_1$
    and $\u_2$, respectively.}  \figlab{comparison1}
\end{figure}

As can be seen, the RMHMC chain is the best in terms of mixing by
comparing the second column (the trace plot) and the third column (the
autocorrelation function ACF). Each RMHMC sample is almost uncorrelated to
the previous ones. The sRMHMC is the second best, but the samples are
strongly correlated compared to those of RMHMC, e.g. one uncorrelated
sample for every $40$. It is interesting to observe that the full RMMALA and sRMMALA have performance in terms of auto-correlation length that is qualitatively similar at least in the first $5000$
samples. This is due to the RMMALA schemes being driven by a single step random walk that cannot exploit fully the curvature information available to the geodesic flows of RMHMC, see rejoinder of \cite{GirolamiCalderhead11}.

 Note that it is not our goal to compare the
behavior of the chains when they converge. Rather we would like to
qualitatively study how fast the chains are well-mixed (mixing time). This is
important for large-scale problems governed by PDEs since
``unpredicted'' mixing time implies a lot of costly waste in PDE solves which
one must avoid. Though RMHMC is expensive in generating a
sample, the cost of generating an uncorrelated/independent sample
seems to be comparable to sRMHMC for this example. In fact, if we
measure the cost in terms of the number of PDE solves, the total number of PDE
solves for RMHMC is $42476480 $ while it is $ 1020002$ for sRMHMC, a
factor of $40$ more expensive. However, the cost in generating an
almost uncorrelated/independent sample is the same since sRMHMC generates one
uncorrelated sample out of $40$ while it is one out of one for RMHMC.

To see how each method distributes the samples we plot one for every
five samples in Figure \figref{comparison1Trajectory}. All methods
seem to explore the high probability density region very
well. This explains why the sample mean and the $95\%$ credibility
region are similar for all methods in the first column of Figure
\figref{comparison1}.

\begin{figure}[h!t!b!]
    \includegraphics[trim=1cm 6.5cm 2cm 6.5cm,clip=true,width=0.97\columnwidth]{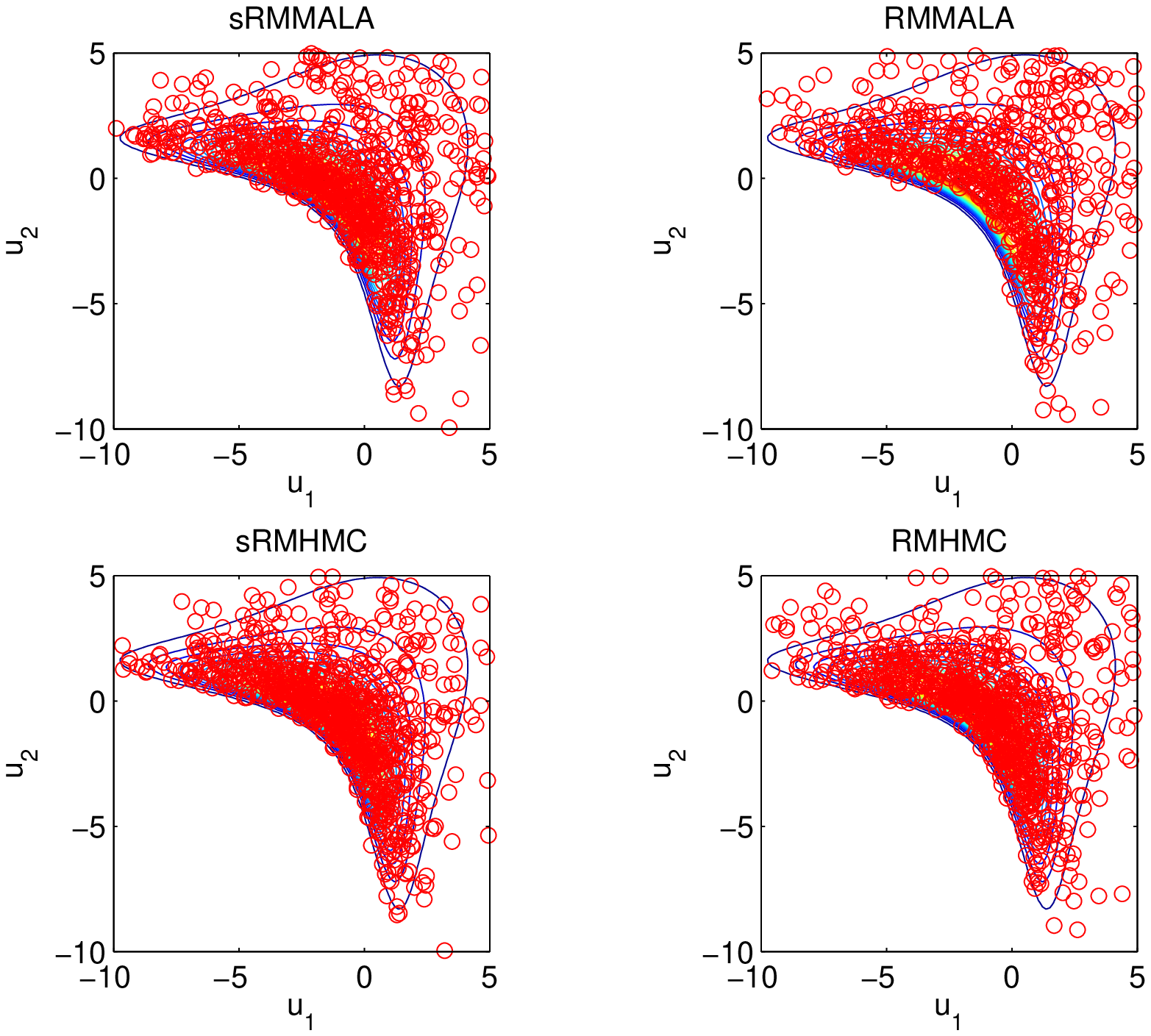}
  \caption{Comparison MCMC trajectories (plot one for each five
    samples) among simRMMALA, RMMALA, simRMHMC, and RMHMC: chains with
    5000 samples, burn-in of $100$, starting at the MAP point. In this
    example, $s = 0.6$, $\alpha = 0.1$, and $\sigma = 0.1$. Time step
    is $\varepsilon = 1$ for simRMMALA and RMMALA, and $\varepsilon =
    0.02$ with the number of time steps $L = 100$ for simRMHMC and
    RMHMC.} \figlab{comparison1Trajectory}
\end{figure}

\begin{figure}[h!tb]
    \includegraphics[trim=1cm 6.5cm 2cm 6.5cm,clip=true,width=0.97\columnwidth]{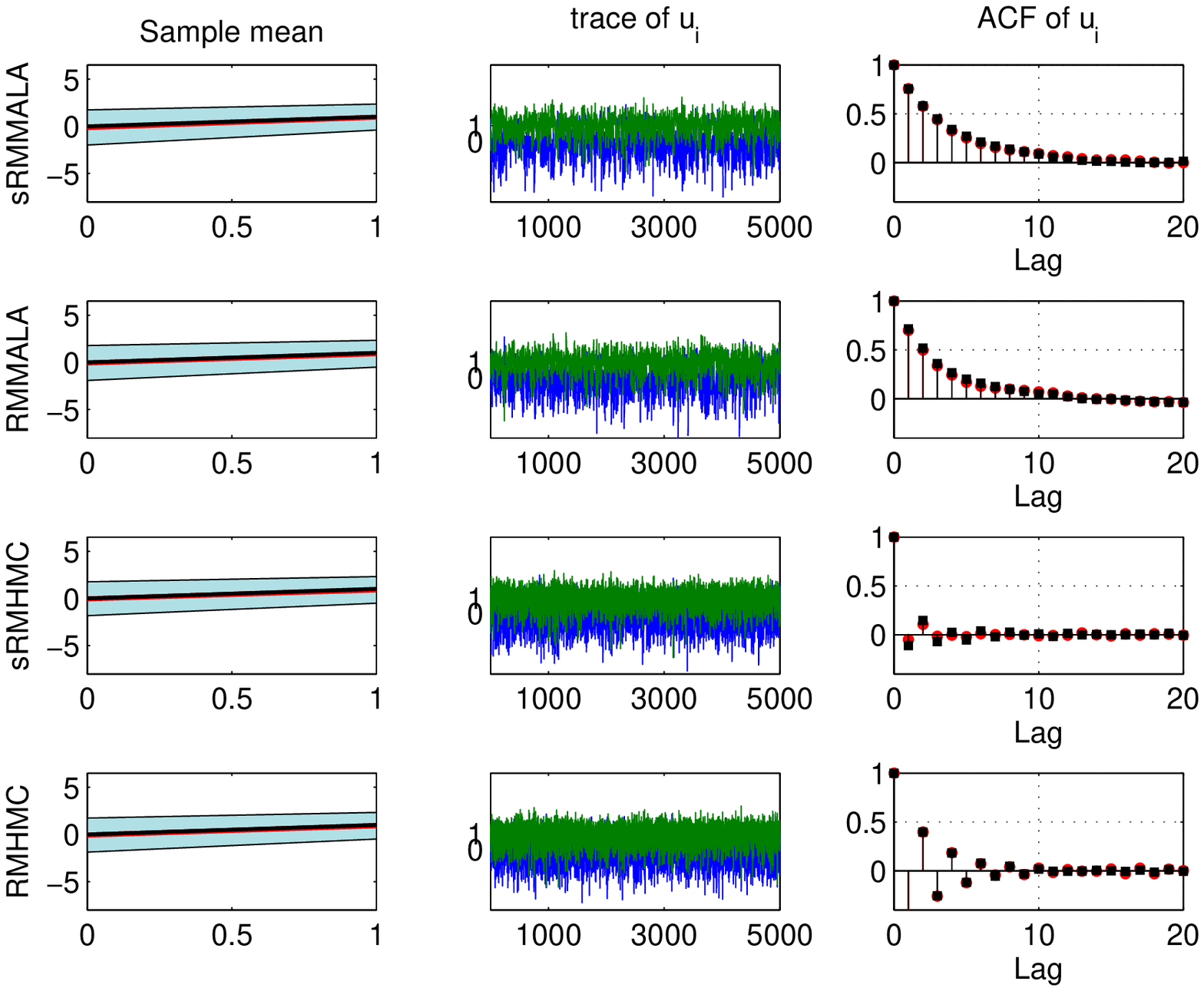}
    \caption{Comparison among simRMMALA, RMMALA, simRMHMC, and RMHMC:
    chains with 5000 samples, burn-in of $100$, starting at the MAP
    point. In this example, $s = 0.6$, $\alpha = 1$, and $\sigma =
    0.01$. Time step is $\varepsilon = 1$ for simRMMALA and RMMALA, and
    $\varepsilon = 0.04$ with the number of time steps $L = 100$ for
    simRMHMC and RMHMC. In the left column: the exact synthetic
    solution is black, the sample mean is red, and the shaded region
    is the $95\%$ credibility region. In the middle column: blue is
    the trace plot for $\u_1$ while green is for $\u_2$. In the right
    column: red and black are the autocorrelation function for $\u_1$
    and $\u_2$, respectively.}
  \figlab{comparison2}
\end{figure}

In the second example we consider the combination $s = 0.6$, $\sigma =
0.01$, and $\alpha = 1$ which leads to the posterior shown in Figure
\figref{posterior2}. For sRMHMC and RMHMC, we take time step
$\varepsilon= 0.04$ with $L = 100$ time steps, while it is $1$ for
both sRMMALA and RMMALA. Again, the acceptance rate is unity for both
sRMHMC and RMHMC while it is $0.65$ for sRMMALA and $0.55$ for RMMALA,
respectively.  The result for four methods is shown in Figure
\figref{comparison2}. As can be seen, this example seems to be easier
than the first one since even though the time step is larger, the
trace plot and the ACF looks better. It is interesting to observe that
sRMHMC is comparable with RMHMC (in fact the ACF seems to be a bit
better) for this example. As a result, RMHMC is more expensive than
sRMHMC for less challenging posterior in Figure
\figref{posterior2}. Here, by less challenging we mean that the
posterior is quite well approximated by a Gaussian at the MAP point,
e.g. the metric tensor is almost constant. This is true for the
posterior in Figure \figref{posterior2} in which the Gaussian prior
contribution is significant, i.e., $\alpha = 1$ instead of $\alpha =
0.1$. Conversely, the posterior is challenging if the metric tensor
changes rapidly. Similar to the first example, one also see that the
sample mean and the $95\%$ credibility region are almost the same for
all methods.



\begin{figure}[h!tb]
    \includegraphics[trim=1cm 6.5cm 2cm 6.5cm,clip=true,width=0.97\columnwidth]{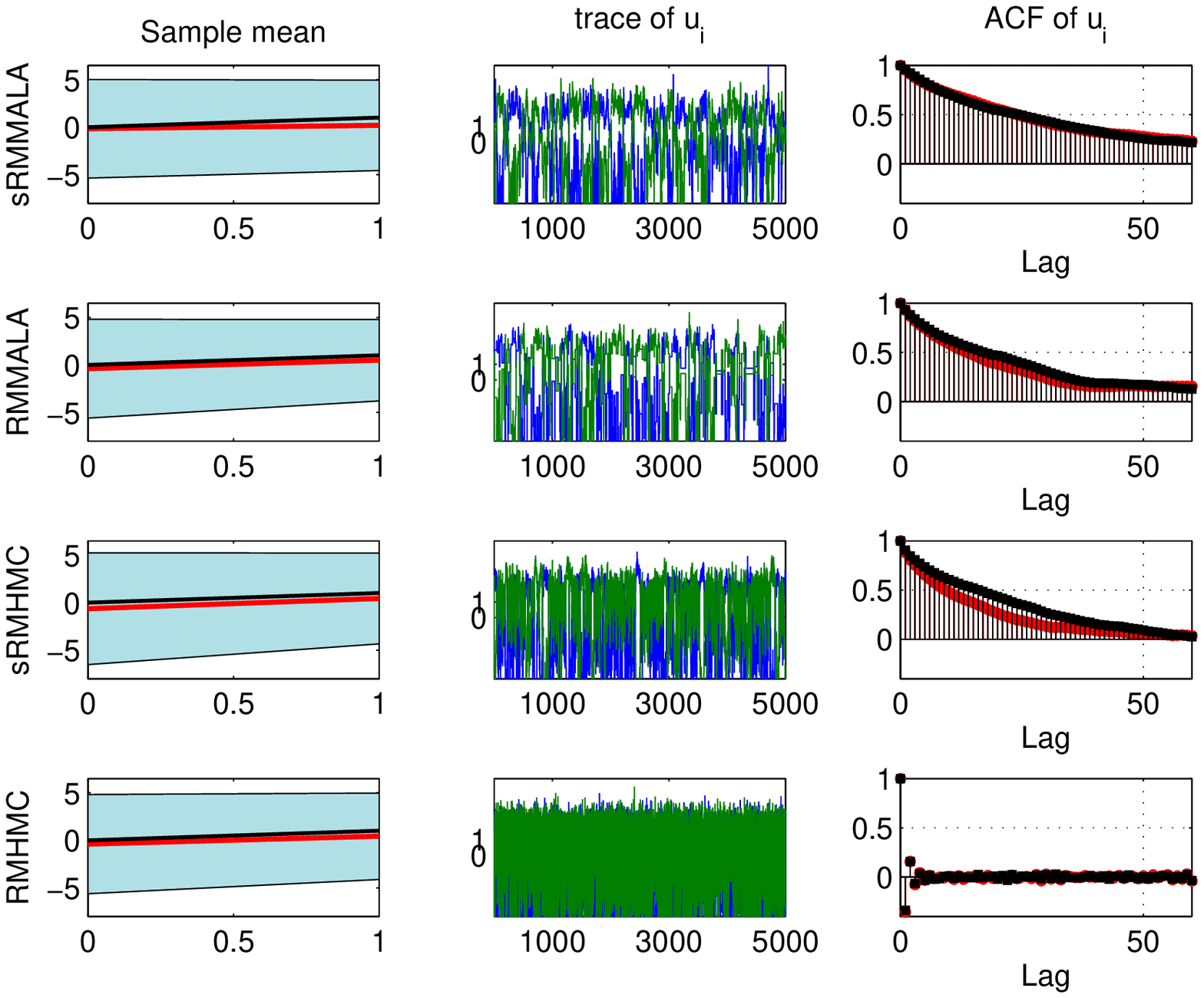}
  \caption{Comparison among simRMMALA, RMMALA, simRMHMC, and RMHMC:
    chains with 5000 samples, burn-in of $100$, starting at the MAP
    point. In this example, $s = 0.6$, $\alpha = 0.1$, and $\sigma =
    0.01$. Time step is $\varepsilon = 0.7$ for simRMMALA and RMMALA, and
    $\varepsilon = 0.02$ with the number of time steps $L = 100$ for
    simRMHMC and RMHMC. In the left column: the exact synthetic
    solution is black, the sample mean is red, and the shaded region
    is the $95\%$ credibility region. In the middle column: blue is
    the trace plot for $\u_1$ while green is for $\u_2$. In the right
    column: red and black are the autocorrelation function for $\u_1$
    and $\u_2$, respectively.}
  \figlab{comparison3}
\end{figure}

In the third example we consider the combination $s = 0.6$, $\sigma =
0.01$, and $\alpha = 0.1$ which leads to a skinny posterior with a
long ridge as shown in Figure \figref{posterior3}. For sRMHMC and
RMHMC, we take time step $\varepsilon= 0.02$ with $L = 100$ time
steps, while it is $1$ for both sRMMALA and RMMALA. Again, the
acceptance rate is unity for both sRMHMC and RMHMC while it is $0.45$
for sRMMALA and RMMALA.  The result for four methods is shown in
Figure \figref{comparison3}. For this example, the RMHMC is more
desirable than sRMHMC since the cost to generate an
uncorrelated/independent sample is smaller for the former than the
latter. The reason is that the total number of PDEs solves for the
former is $40$ times more than the latter, but one out of
very sixty samples is uncorrelated/independent. 


\subsection{Multi-parameter examples}
In this section we choose to discretize $\Omega = \LRs{0,1}$ with
$2^{10} = 1024$ elements, and hence the number of parameters is
$1025$. For all simulations in this section, we choose $s = 0.6$,
$\alpha = 10$, and $\sigma = 0.01$. For synthetic observations, we
take $K = 64$ observations at $x_j = (j-1)/2^6$, $j =
1,\hdots,K$. Clearly, using the full blown RMHMC is out of the
question since it is too expensive to construct the third derivative
tensor and Newton method for each Stomer-Verlet step. For that reason,
the sRMHMC becomes the viable choice.  As studied in Section
\secref{twoparameter}, though sRMHMC loses the ability to efficiently
sample from highly nonlinear posterior surfaces compared to the full
RMHMC it is much less expensive to generate a sample since it does not
require the derivative of the Fisher information matrix. In fact
sRMHMC requires to (approximately) compute the Fisher information at
the MAP point and then uses it as the fixed constant metric tensor
throughout all leap-frog steps for all samples. Clearly, the gradient
\eqnref{Gradient} has to be evaluated at each leap-frog step, but it
can be computed efficiently using the adjoint method presented in
Section \secref{adjoint}.

Nevertheless, constructing the exact Fisher information matrix
requires $2\times 1025$ PDEs solves. This is impractical if the
dimension of the finite element space increases, e.g. by refining the
mesh. Alternatively, due to the compactness of the Hessian of the
prior-preconditioned misfit as discussed in Section \secref{lowrank},
we can use the randomized singular value decomposition (RSVD)
technique \cite{HalkoMartinssonTropp11} to compute its low rank
approximations. Shown in Figure
\figref{eigenSpectrum} are the first 35 dominant eigenvalues of the
Fisher information matrix and its prior-preconditioned counterpart. We
also plot 20 approximate eigenvalues of the prior-preconditioned
Fisher information matrix obtained from the RSVD method. As can be
seen, the eigen spectrum of the prior-preconditioned Fisher
information matrix decays faster than that of the original one. This
is not surprising since the prior-preconditioned Fisher operator is a
composition of the prior covariance, a compact operator, and the
Fisher information operator, also a compact operator. The power of the
RSVD is clearly demonstrated as the RSVD result for the first 20
eigenvalues is very accurate.

\begin{figure}[h!tb]
    \includegraphics[trim=1cm 6.5cm 2cm 6.5cm,clip=true,width=0.97\columnwidth]{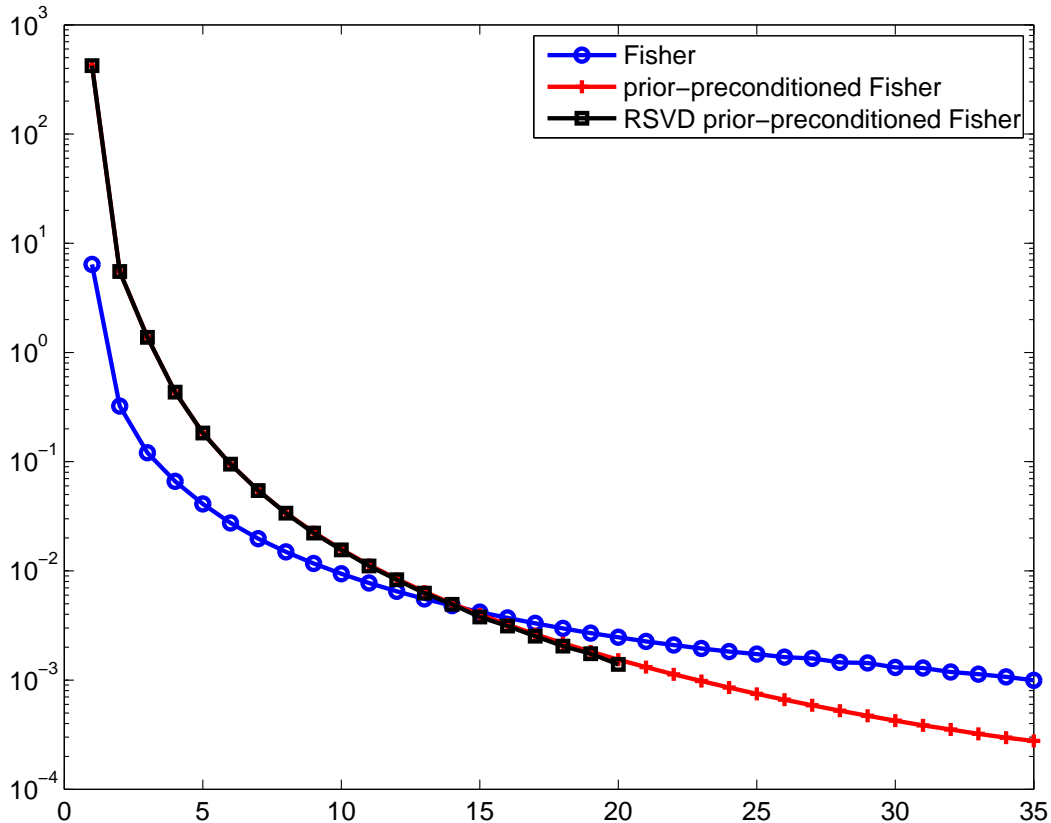}
  \caption{The eigen spectrum of the Fisher information matrix, the
    prior-preconditioned Fisher matrix, and the first 20 eigenvalues
    approximated using RSVD. Here, $s = 0.6$, $\alpha = 10$, and $\sigma
    = 0.01$.}  \figlab{eigenSpectrum}
\end{figure}

\begin{figure}[h!tb]
    \includegraphics[trim=1cm 6.5cm 2cm 6.5cm,clip=true,height = 0.7\columnwidth, width=1\columnwidth]{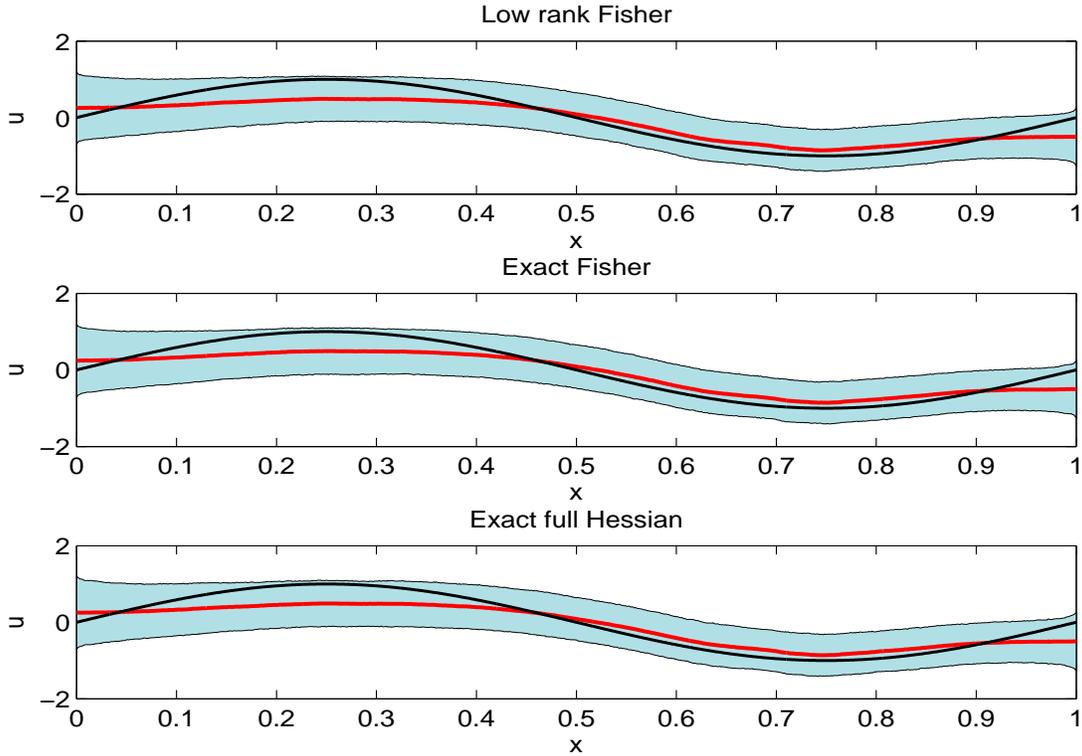}
  \caption{MCMC results of three sRMHMC method with i) the low rank
    Gauss-Newton Hessian, ii) the exact Gauss-Newton Hessian, and iii)
    the full Hessian. In the figure are the empirical mean (red line),
    the exact distributed parameter used to generate the observation
    (black line), and $95\%$ credibility (shaded region).}
  \figlab{UQmultiparameter}
\end{figure}

Next, we perform the sRMHMC method using three different constant
metric tensors: i) the low rank Gauss-Newton Hessian, ii) the exact
Gauss-Newton Hessian, and iii) the full Hessian. For
each case, we start the Markov chain at the MAP point and compute
$5100$ samples, the first $100$ of which is then discarded as burn-ins. The
empirical mean (red line), the exact distributed parameter used to
generate the observation (black line), and $95\%$ credibility region
are shown in Figure \figref{UQmultiparameter}. As can be seen, the
results from the three methods are indistinguishable. The first sRMHMC
is the most appealing since it requires $2 \times 20 = 40$ PDE solves
to construct the low rank Fisher information while the others need
$2\times 1025$ PDE solves. For large-scale problems with computationally
expensive PDE solves, the first approach is the method of choice.

To further compare the three methods we record the trace plot of the
first two ($1$ and $2$) and the last two ($1024$ and $1025$) parameters in
Figure \figref{UQmultiparameterTrace}. As can be observed, the chains
from the three methods seem to be well-mixed and it is hard to see the
difference among them. We also plot the autocorrelation function for
these four parameters. Again, results for the three sRMHMC methods are
almost identical, namely, they generate almost uncorrelated samples. We
therefore conclude that low rank approach is the least computational
extensive, yet it maintains the attractive features of the original
RMHMC. As such, it is the most suitable method for large-scale
Bayesian inverse problems with costly PDE solves.

\begin{figure}[h!t!b!]
    \includegraphics[trim=1cm 6.5cm 2cm 6.5cm,clip=true,height = 1.\columnwidth, width=1\columnwidth]{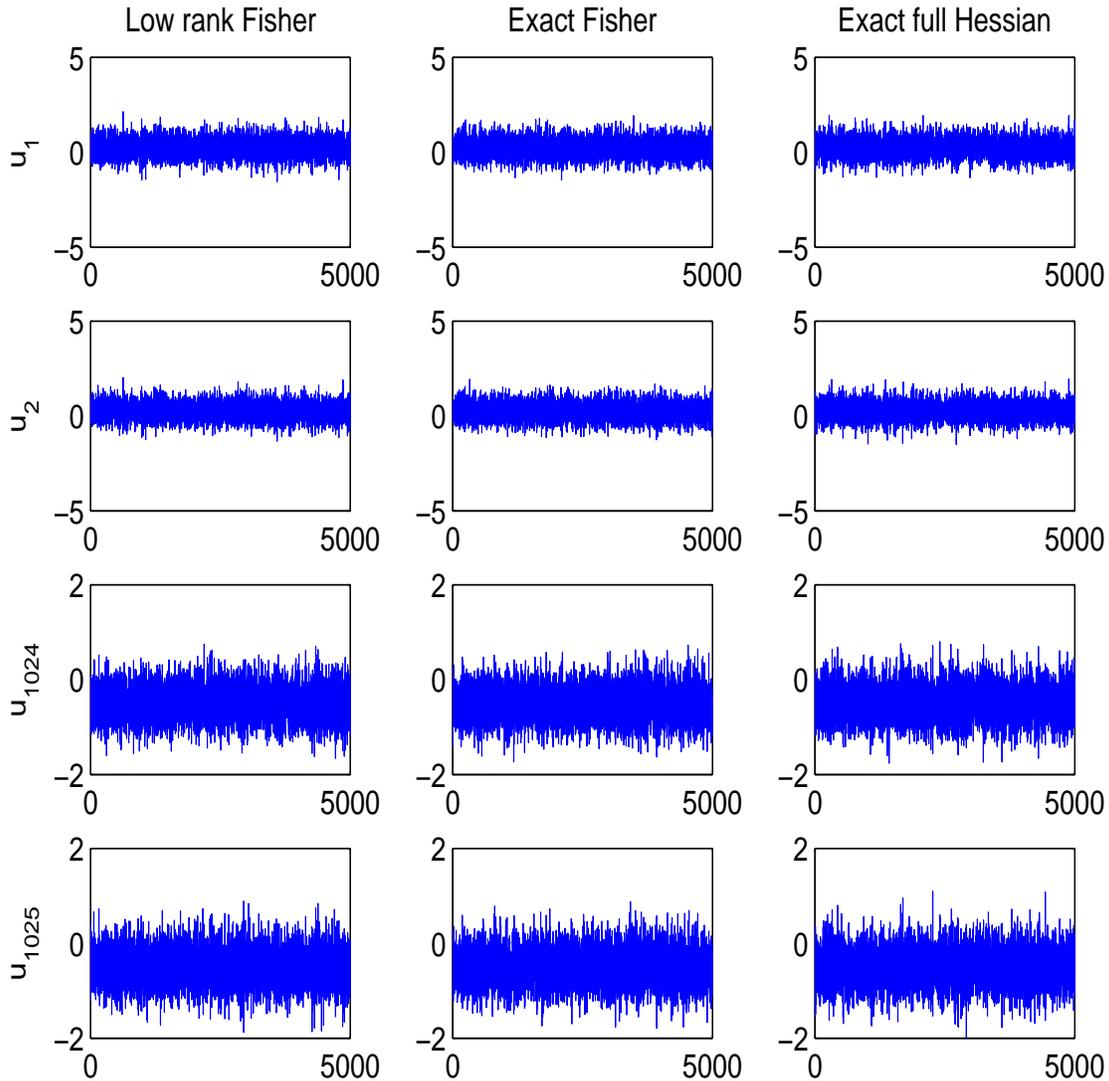}
  \caption{MCMC results of three sRMHMC method with i) the low rank
    Gauss-Newton Hessian (left column), ii) the exact Gauss-Newton
    Hessian (middle column), and iii) the full Hessian at the MAP
    point (right column). In the figure are the trace plot of the
    first two ($1$ and $2$) and last two ($1024$ and $1025$) parameters}
  \figlab{UQmultiparameterTrace}
\end{figure}

\begin{figure}[h!t!b!]
    \includegraphics[trim=1cm 6.5cm 2cm 6.5cm,clip=true,height = 1.\columnwidth, width=1\columnwidth]{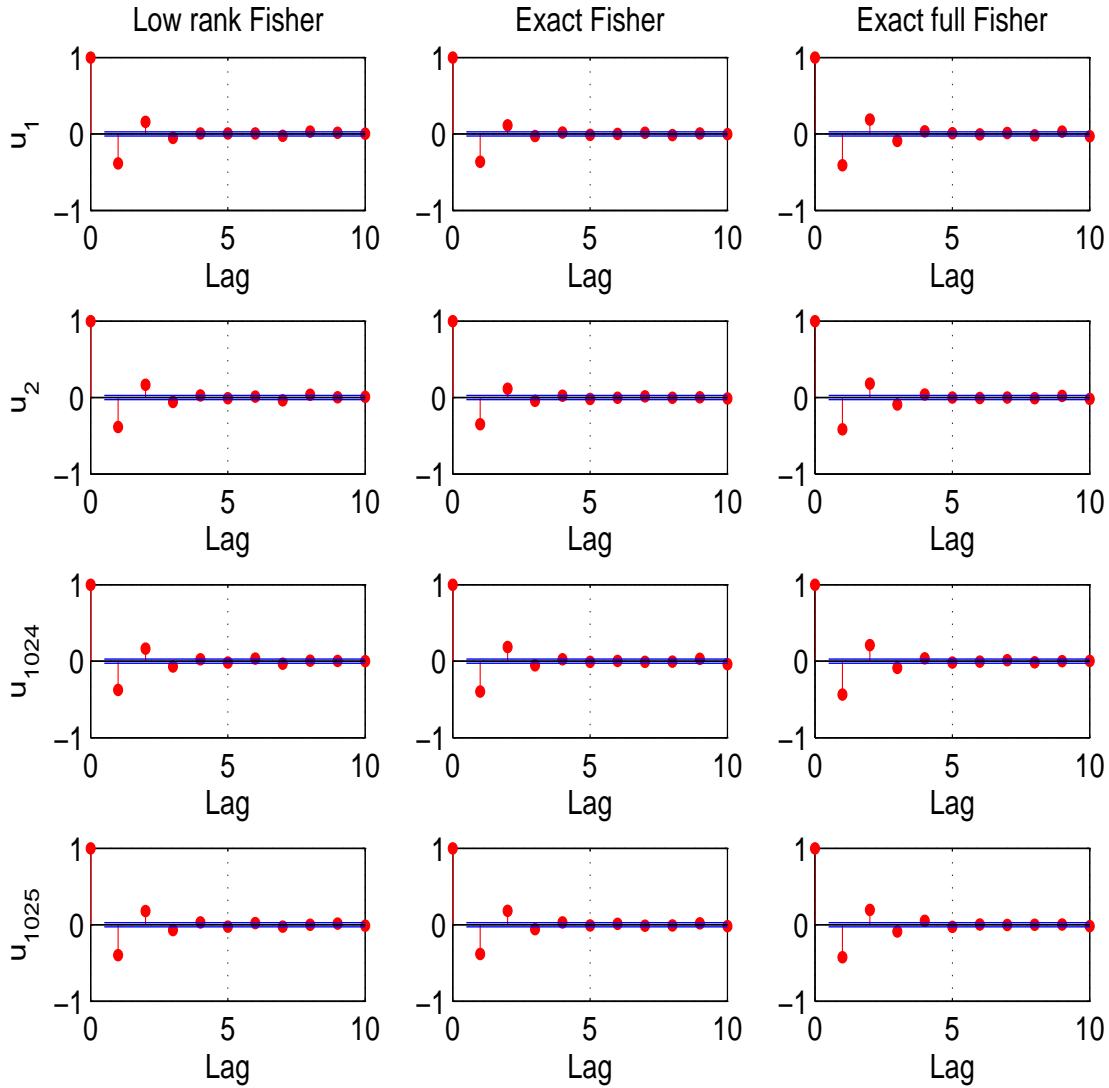}
  \caption{MCMC results of three sRMHMC method with i) the low rank
    Gauss-Newton Hessian (left column), ii) the exact Gauss-Newton
    Hessian (middle column), and iii) the full Hessian at the MAP
    point (right column). In the figure are the autocorrelation function plot of the
    first two ($1$ and $2$) and last two ($1024$ and $1025$) parameters}
  \figlab{UQmultiparameterCorrelation}
\end{figure}



\section{Conclusions and future work}
\seclab{conclusions} We have proposed the adoption of a computationally  inexpensive Riemann manifold
Hamiltonian Monte method to explore the posterior of large-scale
Bayesian inverse problems governed by PDEs in a highly efficient manner. We first adopt an infinite
dimensional Bayesian framework to guarantee that the inverse
formulation is well-defined. In particular, we postulate a Gaussian
prior measure on the parameter space and assume regularity for the
likelihood. This leads to a well-defined posterior distribution. Then,
we discretize the posterior using the standard finite element method
and a matrix transfer technique, and apply the RMHMC method on the
resulting discretized posterior in finite dimensional parameter
space. We present an adjoint technique to efficiently compute the
gradient, the Hessian, and the third derivative of the potential
function that are required in the RMHMC context. 
This is at the expense of
solving a few extra PDEs: one for the gradient, two for a
Hessian-vector product, and four for the product of third order
derivative with a matrix.
For large-scale
problems, repeatedly computing the action of the Hessian and third order derivative
is too computationally expensive and this motivates us to design a
simplified RMHMC in which the Fisher information matrix is computed
once at the MAP point. We further reduce the effort by constructing
low rank approximation of the Fisher information using a randomized singular value decomposition technique. The effectiveness of
the proposed approach is demonstrated on a number of numerical results
up to $1025$ parameters in which the computational gain is about two
orders of magnitude while maintaining the quality of the original
RMHMC method in generating (almost) uncorrelated/independent samples.

For more challenging inverse problems with significant change of
metric tensor across the parameter space, we expect that sRMHMC with
constant metric tensor is inefficient. In that case, RMHMC seems to be a
better option, but it is too computational extensive for large-scale
problems. Ongoing work is to explore approximation of the RMHMC
methods in which we approximate the trace and the third derivative in
\eqnref{Hamilton} using adjoint and randomized techniques.

\section*{References}
\bibliographystyle{unsrt}

\end{document}